\documentclass[11pt]{article}

\usepackage[utf8]{inputenc}
\usepackage[T1]{fontenc}
\usepackage{lmodern}
\usepackage[margin=1in]{geometry}
\usepackage{amsmath,amssymb,amsthm}
\usepackage{graphicx}
\usepackage{authblk}
\usepackage{tikz}
\usetikzlibrary{arrows.meta, positioning, fit, backgrounds}
\usepackage{multirow}
\usepackage{url}

\newcommand{\R}{\mathrm{I\hspace{-0.5ex}R}}

\newcommand{\G}{\mathcal{G}}

\usepackage{float}

\title{High-performance computing for the BGK model of the Boltzmann equation
	with a meshfree Arbitrary Lagrangian--Eulerian method}

\author[1]{Panchatcharam Mariappan}
\author[2]{Klaas Willems}
\author[3]{Gangadhara Boregowda}
\author[2]{Sudarshan Tiwari}
\author[2]{Axel Klar}

\affil[1]{Department of Mathematics and Statistics, Indian Institute of
	Technology Tirupati, Tirupati, Andhra Pradesh, India
	(\texttt{panch.m@iittp.ac.in})}
\affil[2]{Department of Mathematics, RPTU Kaiserslautern-Landau,
	Erwin-Schr\"odinger-Stra{\ss}e, 67663 Kaiserslautern, Germany}
\affil[3]{Department of Mathematics, SRM Institute of Science and Technology,
	Kattankulathur, Tamil Nadu, India}
\date{}

\begin{document}
	
	\maketitle
	
\begin{abstract}
	In this paper, we present high-performance computing for the BGK model of the
	Boltzmann equation with a mesh-free method. For the numerical simulation of the
	BGK equation, we use an Arbitrary-Lagrangian-Eulerian (ALE) method developed in
	previous work, where the approximation of spatial derivatives and the
	reconstruction of a function are based on the weighted least squares method. A
	Graphics Processing Unit (GPU) is used to accelerate the code. The result is
	compared with sequential and parallel versions of the CPU code. Two- and
	three-dimensional driven cavity problems are solved, where a speed-up of
	several orders of magnitude is obtained compared to a sequential CPU
	simulation. The transport step is implemented as a node-parallel CUDA kernel
	for both two- and three-dimensional physical space: one GPU thread advances one
	particle-velocity-node pair by tracing its characteristic foot point and
	reconstructing the distribution there with a register-resident weighted least
	squares solve. For the three-dimensional case, no velocity reduction is
	applicable: a single distribution function $f$ is evolved on the full $N_v^3$
	velocity grid, so the transport in every coordinate direction follows the
	resolved velocity node. The device memory footprint is kept at exactly twice
	the distribution storage by a structure-of-arrays layout with a single
	transport buffer; particle management runs on the host, so no device-side
	scratch copies of the cloud are required.
\end{abstract}
	
\section{Introduction}

During the last two decades, simulation of Micro Electro Mechanical Systems (MEMS) has attracted many researchers; see \cite{KBA, RF, DM, TKHD, STKH, STK, BTSKH, TA, TA1}. This is mostly micro- and nanofluidics, where the mean free path is often of the order of or larger than the characteristic length of the geometry. That means the Knudsen number, i.e., the ratio of the mean free path to the characteristic length, is of the order of 1 or larger.  In this case, one has to solve rarefied gas flows modelled by the Boltzmann equation. Usually, these flows have low Mach numbers; therefore, stochastic methods like Direct Simulation Monte Carlo (DSMC)~ \cite{Bird} are not the optimal choice, since statistical noise dominates the flow quantities. Moreover, when a rigid body moves through the gas, the domain changes in time and one encounters unsteady flow problems such that averages over long runs cannot be taken. Instead, one has to perform many independent runs and average them to obtain the smooth solutions, which is time-consuming.
Although some attempts have been made to reduce the statistical noise of DSMC-type methods, see, for example, \cite{DDP}, or to adopt efficient solvers for the Boltzmann equation, such as those based on the Fourier-spectral method (see, for example, the review paper \cite{dimarco-acta2014}), many works rather employ deterministic approaches for simplified models of the Boltzmann equation, like the Bhatnagar-Gross-Krook (BGK) model; see \cite{RF, DM, TA, TKR, TKR2}.

This paper is based on an earlier paper \cite{TKR3} describing a mesh-free Arbitrary Lagrangian-Eulerian (ALE) method for the BGK equation. There, we have solved two-dimensional problems with moving objects inside a rarefied flow. The computational costs for such simulations are high due to the high dimensionality of the problem. This is particularly true for 3D problems where, in addition to the three spatial dimensions, a further 3 velocity dimensions have to be considered for the BGK equations.
Thus, reasonably accurate numerical simulations require a suitable parallelization procedure. In the present paper, such a procedure using Graphics Processing Unit (GPU) computations and Message Passing Interface (MPI) is presented and numerically investigated.

In this paper, we exploit the GPU, which is well-known for massively parallel computing. Researchers have developed and demonstrated the advantage of GPUs for kinetic models using various schemes such as DSMC \cite{Frezzotti2, Frezzotti3, Frezzotti1}, the Lattice Boltzmann method \cite{Vuyst1, Vuyst2, Tadeusz, Malkov, Andrea}, discrete unified gas kinetic scheme \cite{Peiyao}, direct numerical simulation \cite{Yuhang}, and a semi-Lagrangian approach \cite{DimarcoIII}.

A key contribution of the present work is the node-parallel GPU implementation of the meshfree transport step in both 2D and 3D physical space.  Parallelism is over particle--velocity-node pairs ($N\cdot N_v^2$ threads in 2D, $N\cdot N_v^3$ in 3D) rather than over particles alone, so the abundant fine-grained work of the kinetic velocity grid saturates the GPU.  Each thread solves its weighted least squares reconstruction entirely in registers with a compile-time system size, and the global conservation correction uses an atomic-free block-reduction scheme -- on the A100 this is decisive, as a single-address floating-point atomic reduction serialises and was measured to cost three orders of magnitude more than the block-parallel equivalent.  Velocity-grid consistency requires a sufficiently wide cutoff: we show that a cutoff of at least $5\,\sigma$ (where $\sigma=\sqrt{RT}$ is the thermal speed) is necessary to avoid a discretisation-induced temperature mismatch between the bulk equilibrium and the wall boundary condition that otherwise seeds a slowly growing instability.

The paper is organized as follows. In section \ref{sec:model}, we present the BGK model for the Boltzmann equation including the Chu reduction procedure used for 2D computations. In Section \ref{sec:methods}, we present the numerical scheme for the BGK model, in particular the spatial and temporal discretization using a first-order discretization method. In Section \ref{sec:HPC}, we present the GPU architecture and the implementation of the numerical scheme on the GPU. In Section \ref{sec:numerics}, we present numerical results in two- and three-dimensional cases. Finally, in Section \ref{sec:conclusion}, some conclusions and an outlook are presented.

\section{The BGK model for rarefied gas dynamics}
\label{sec:model}

We consider the BGK  model of the Boltzmann equation for rarefied gas dynamics, where the collision term is modeled by a relaxation of the distribution function $f(t,{\bf x},{\bf v})$ to the Maxwellian equilibrium distribution. The evolution equation for the distribution function $f(t,{\bf x},{\bf v})$  is given by the following initial boundary value problem
\begin{equation}
	\frac{\partial f}{\partial t} + {\bf v}\cdot\nabla_{\bf x} f = \frac{1}{\tau}(M - f),
	\label{bgk_eqn}
\end{equation}
with $f(0,{\bf x}, {\bf v}) = f_0({\bf x}, {\bf v}), \; t \ge 0, {\bf x} \in \Omega \subset \mathbb R^{d_x}, (d_x=1,2,3), \;  {\bf v}  \in \mathbb R^{d_v}, (d_v=1,2,3)$  and suitable initial and boundary conditions described in the next section.
Here $\tau$ is the relaxation time and $M$ is the local Maxwellian given by
\begin{equation}
	M = \frac{\rho}{(2\pi R T)^{d_v/2}} \exp \left(-\frac{| {\bf v} -  {\bf U}|^2}{2RT}\right),
	\label{maxwellian}
\end{equation}
where the parameters $\rho (t,{\bf x}) \in \R, {\bf U} (t,{\bf x}) \in \R^{d_v}, T (t,{\bf x})\in \R$ are the macroscopic quantities density, mean velocity, and temperature, respectively.  $R$ is the specific gas constant of the gas.  $\rho, {\bf U} ,T $ are computed from $f$ as  follows.
Let the moments of $f$ be  defined by
\begin{equation}
	(\rho, \rho {\bf U}, E) = \int_{\mathbb R^{3}} \phi({\bf v}) f(t,  {\bf x},  {\bf v}) d {\bf v}.
	\label{moments}
\end{equation}
where  $ { \phi } ({\bf v})=\left (1,  {\bf v} ,\frac{| {\bf v} |^2}{2} \right )$ denotes  the vector of collision invariants.
$E$ is the total energy density, which is related to the temperature through the internal energy
\begin{equation}
	e(t, {\bf x}) = \frac{3}{2}R T, \quad \quad \rho e = E - \frac{1}{2}\rho |{\bf U}|^2.
	\label{internal_energy}
\end{equation}
The relaxation time $\tau=\tau(x,t)$ and the mean free path $\lambda$ are related according to  \cite{CC}
\begin{equation}
	\tau = \frac{4 \lambda}{\pi \bar C},
	\label{tau}
\end{equation}
where $\bar C = \sqrt{\frac{8RT}{\pi}}$ and the mean free path is given by

\begin{equation}
	\lambda = \frac{k_b}{\sqrt{2}\pi\rho R d^2},
	\label{lambda}
\end{equation}
where $k_b$ is the Boltzmann constant and $d$ is the diameter of the gas molecules.

\subsection{Chu-reduction}
To solve two-dimensional flow problems with $d_x = 2$, one might consider a two-dimensional velocity space $d_v = 2$ mathematically, since the consideration of three-dimensional velocity space requires unnecessarily high memory and computational time.
However,   it is physically correct to consider in these situations still three-dimensional velocity space.
In these cases, for the BGK model,
the  3D velocity space can be reduced as suggested by Chu \cite{Chu}. This reduction yields a considerable reduction in memory allocations and computational time.
For example, in a  physically two-dimensional situation, in which all variables depend on ${\bf x} \in \R^2$ and $t$,  the
velocity space is reduced from three dimensions to two dimensions, defining the following reduced distributions \cite{GR}. Considering  ${\bf v}  = (v_1,v_2,v_3) \in \R^3$ we define
\begin{eqnarray}
	g_1(t,x,v_1, v_2) &=& \int_{\mathbb R} f(t,x,v_1,v_2,v_3)dv_3
	\label{reduced_g1}
	\\
	g_2(t,x,v_1, v_2) &= &\int_{\mathbb R} v_3^2 f(t,x,v_1,v_2,v_3)dv_3.
	\label{reduced_g2}
\end{eqnarray}
Multiplying (\ref{bgk_eqn}) by $1$ and $v_3^2$ and integrating with respect to $v_3 \in \mathbb R$, we obtain the following
system of two equations denoting for  simplicity $(v_1, v_2) = {\bf v}$ in the reduced equations.
\begin{eqnarray}
	\frac{\partial g_1}{\partial t} + {\bf v}\cdot \nabla_{\bf x} g_1 &=& \frac{1}{\tau}(G_1 - g_1)
	\label{eqduced_bgk1}
	\\
	\frac{\partial g_2}{\partial t} + {\bf v}\cdot \nabla_{\bf x} g_2 &=& \frac{1}{\tau}(G_2 - g_2),
	\label{eqduced_bgk2}
\end{eqnarray}
where
\begin{eqnarray}
	G_1 &= & \int_{\mathbb R} M dv_3 = \frac{\rho\pi}{2RT} \exp \left(-\frac{| {\bf v} -  {\bf U}|^2}{2RT}\right)
	\label{Maxw_G1}
	\\
	G_2 &= & \int_{\mathbb R} v_3^2 M dv_3 = (RT) G_1
	\label{Maxw_G2}
\end{eqnarray}
with ${\bf U}= (U_1,U_2) \in \R^2 $.
Assuming that it  is a local equilibrium, the initial distributions are defined via the parameters $\rho_0, {\bf U}_0, T_0$ and are given as
\begin{eqnarray}
	g_1(0,{\bf x},{\bf v}) =  \frac{\rho_0\pi}{2RT_0} \exp \left(-\frac{| {\bf v} -  {\bf U}_0|^2}{2RT_0}\right), \quad
	g_2(0,{\bf x}, {\bf v}) = (RT_0) g_1(0,{\bf x}, {\bf v}).
\end{eqnarray}
The macroscopic quantities are given through the reduced distributions as
\begin{equation}
	\rho = \int_{\mathbb R^2} g_1d{\bf v}, \; \rho {\bf U} = \int_{\mathbb R^2}{\bf v} g_1 d{\bf v}, \; 3\rho RT = \int_{\mathbb R^2}|{\bf v}-{\bf U}|^2 g_1 d{\bf v} + \int_{\mathbb R^2} g_2 d{\bf v}.
	\label{macroeq}
\end{equation}
%


In a three-dimensional simulation, all
three velocity components are physical degrees of freedom, so no reduction
is applicable and a single distribution function $f(t,{\bf x},{\bf v})$ is
evolved on the full three-dimensional velocity grid. All macroscopic
quantities are then recovered from this one distribution,
\begin{equation}
	(\rho,\,\rho {\bf U},\, E)(t,{\bf x})
	= \int_{\R^3} \bigl(1,\,{\bf v},\,\tfrac{1}{2}|{\bf v}|^2\bigr)\,
	f(t,{\bf x},{\bf v})\, d{\bf v}, \qquad
	T = \frac{1}{3R}\Bigl(\frac{2E}{\rho} - |{\bf U}|^2\Bigr)
\end{equation}
and the relaxation is toward the full Maxwellian
\begin{equation}
	M[f] \;=\; \frac{\rho}{(2\pi R T)^{3/2}}
	\exp \Bigl(-\frac{|{\bf v}-{\bf U}|^2}{2RT}\Bigr).
\end{equation}
The memory per grid point grows as $N_v^{3}$, so the
three-dimensional computations use moderate $N_v$. The transport step
follows the true characteristic in all three directions; in particular, a
wall moving normal to the third coordinate transfers momentum correctly,
which the reduced pair $(g_1,g_2)$ cannot represent.

\section{Numerical schemes}\label{sec:methods}

We solve the original equation (\ref{bgk_eqn}) in three-dimensional physical space and the reduced system of equations (\ref{eqduced_bgk1} - \ref{eqduced_bgk2}) in two-dimensional physical space by the ALE method described below. We use a time splitting, where the
advection step is solved explicitly, and the relaxation part is solved implicitly, and a meshfree particle method is used to solve this system of equations while approximating the spatial derivatives.
Here, by particle, we actually mean grid points moving with the mean velocity ${\bf U}$ of the gas. The spatial derivatives
of the distribution function at an arbitrary particle position are approximated using values at the point cloud surrounding the particle and a weighted least squares method.


\subsection{ALE formulation}

The semi-Lagrangian method is based on fixed grids. For flows with moving boundaries, meshfree methods with moving grids might be more suitable than fixed grids; see \cite{TKR2, TKR3}. For such a method, one rewrites the equations as follows.

\subsubsection{ALE formulation for original model}
We rewrite  the equations (\ref{bgk_eqn}) in Lagrangian form as
\begin{eqnarray}
	\frac{d{\bf x}}{dt} &=& {\bf U} \label{originalALE1}
	\\
	\frac{d f}{d t} &=& -({\bf v} - {\bf U}) \cdot \nabla_{\bf x} f+ \frac{1}{\tau}(M-f)
	\label{originalALE2}
\end{eqnarray}
where $\frac{d }{dt} = \frac{\partial }{\partial t} + {\bf U}\cdot \nabla_x$. The first equation describes motion with the macroscopic
mean velocity ${\bf U}$ of the gas determined by (\ref{macroeq}). The second equation includes the remaining advection with the difference between microscopic and macroscopic velocity.

\subsubsection{ALE for reduced  model}

In this case  the equations (\ref{reduced_g1} - \ref{reduced_g2}) are reformulated in  Lagrangian form as
\begin{eqnarray}
	\frac{d{\bf x}}{dt} &=& {\bf U} \label{ALE1}
	\\
	\frac{d g_1}{d t} &=& -({\bf v} - {\bf U}) \cdot \nabla_{\bf x}  g_1+ \frac{1}{\tau}(G_1 - g_1)
	\label{ALE2}
	\\
	\frac{d g_2}{d t} &=& - ({\bf v} - {\bf U}) \cdot \nabla_{\bf x} g_2 + \frac{1}{\tau}(G_2 - g_2).
	\label{ALE3}
\end{eqnarray}

\subsection{Time discretization}

The above equations are discretized with the following scheme.

\subsubsection{First order time splitting scheme for original model}
\label{1st_order_step}
Time is discretized as  $t^n = n \Delta t, n= 0, 1, \cdots, N_t$. We denote the numerical approximation of $f$ at $t_n$ by $f^n = f(t^n, {\bf x}, {\bf v})$.
We use a time splitting scheme for equation (\ref{originalALE2}), where the advection term is solved explicitly, and the collision term is solved implicitly.
In the first step of the splitting scheme, we obtain the intermediate distribution $\tilde{f}^n$ by solving
\begin{eqnarray}
	\tilde{f}^n &=& f^n - \Delta t ({\bf v} - {\bf U}^n) \cdot \nabla_{\bf x} f^n .
	\label{explicit}
\end{eqnarray}
In the second step, we obtain the new distribution by solving
\begin{eqnarray}
	f^{n+1} &=& \tilde{f}^n+ \frac{\Delta t}{\tau}(M^{n+1} - f^{n+1} )
	\label{implicitfull_1}
\end{eqnarray}
and the new positions of the grids are updated by
\begin{equation}
	{\bf x}^{n+1} = {\bf x}^n + \Delta t {\bf U}^{n+1}.
\end{equation}

In the first step, we have to approximate the spatial derivatives of $f$ at every grid point. This is described in the following section.

Following \cite{GRS, GR, XRQ}
we  obtain $f^{n+1}$ in the second step by  first  determining  the parameters $\rho^{n+1}, {\bf U}^{n+1}$ and $T^{n+1}$ for $M^{n+1}$.
Multiplying (\ref{implicitfull_1}) by $1$, ${\bf v}$ and $|{\bf v} - {\bf U} |^2$ and  integrating  over velocity space,  we get
\begin{equation}
	\rho^{n+1} = \int_{\mathbb R^3} \tilde{f}^n d{\bf v}, \quad (\rho {\bf U})^{n+1} = \int_{\mathbb R^3} v\tilde{f}^n d{\bf v}, \quad
	3\rho RT^{n+1} = \int_{\mathbb R^3} |{\bf v} - {\bf U}|^2 \tilde{f}^n d{\bf v},
	\label{fullmoment1a}
\end{equation}
where we have used the conservation of mass, momentum, and energy of the original BGK model.

Now, the parameters $\rho^{n+1}, {\bf U}^{n+1}$ and $T^{n+1}$ of $M^{n+1}$  are given in terms of $\tilde{f}$  from (\ref{fullmoment1a}) since $\rho, U$ and $T$ of $f$ and $M$ are same. Hence the implicit step (\ref{implicitfull_1}) can be rewritten as
\begin{eqnarray}
	f^{n+1} = \frac{\tau \tilde{f}^n + \Delta t M^{n+1}}{\tau + \Delta t}.
\end{eqnarray}

\subsubsection{Time splitting scheme for the reduced model}
\label{1st_order_step_red}

We use a time splitting scheme again. In the first step we obtain the intermediate distributions $\tilde{g}^n_1$  and $\tilde{g}^n_2$ by solving
for ${\bf v} \in {\mathbb R^2}$ and ${\bf U} \in {\mathbb R^2}$
\begin{eqnarray}
	\tilde{g}_1^n &=& g_1^n - \Delta t ({\bf v} - {\bf U}^n) \cdot \nabla_{\bf x}   g_1^n
	\\
	\tilde{g}_2^n &=& g_2^n - \Delta t ({\bf v} - {\bf U}^n) \cdot \nabla_{\bf x} g_2^n .
\end{eqnarray}
In the second step, we obtain the new distributions by solving
\begin{eqnarray}
	g_1^{n+1} &=& \tilde{g}_1^n + \frac{\Delta t}{\tau}(G_1^{n+1} - g_1^{n+1} )
	\label{implicit_1}
	\\
	g_2^{n+1} &=& \tilde{g}_2^n + \frac{\Delta t}{\tau}(G_2^{n+1} - g_2^{n+1} )
	\label{implicit_2}
\end{eqnarray}
and the new positions of the grids are updated by
\begin{equation}
	{\bf x}^{n+1} = {\bf x}^n + \Delta t {\bf U}^{n+1}.
\end{equation}

For the second step,  we have to determine first the parameters $\rho^{n+1}, U^{n+1}$ and $T^{n+1}$ for $G_1^{n+1}$ and $G_2^{n+1}$.
Multiplying (\ref{implicit_1}) by $1$ and $v$ and integrating with respect to $v$ over $\mathbb R^2$ we get
\begin{equation}
	\rho^{n+1} = \int_{\mathbb R^2} \tilde{g}_1^n d{\bf v}, \quad (\rho U)^{n+1} = \int_{\mathbb R^2} v\tilde{g}_1^n d{\bf v},
	\label{moment1a}
\end{equation}
where we have used the conservation of mass and momentum of the original BGK model.  In order to compute $T^{n+1}$
we note that the following identity is valid
\begin{equation}
	\int_{\mathbb R^2} |{\bf v} -{\bf U}|^2(G_1-g_1)d{\bf v} + \int_{\mathbb R^2} (G_2 - g_2) d{\bf v} = 0.
	\label{id3}
\end{equation}

Multiplying the equation (\ref{implicit_1}) by $|{\bf v} - {\bf U}|^2$ and integrate with respect to ${\bf v}$ over $\mathbb R^2$ we get
\begin{equation}
	\int_{\mathbb R^2} |{\bf v} - {\bf U})^2 g_1^{n+1} d{\bf v} = \int_{\mathbb R^2} |{\bf v} - {\bf U}|^2 \tilde{g}_1^n d{\bf v} + \frac{\Delta t}{\tau} \int_{\mathbb R^2} |{\bf v} - {\bf U}|^2(G_1^{n+1} - g_1^{n+1} ) d{\bf v}.
	\label{id4}
\end{equation}

Next,  integrate both sides of (\ref{implicit_2}) with respect to ${\bf v}$ over $\mathbb R^2$ we get
\begin{equation}
	\int_{\mathbb R^2} g_2^{n+1}d{\bf v} = \int_{\mathbb R^2}\tilde{g}_2^n d{\bf v} + \frac{\Delta t}{\tau} \int_{\mathbb 3R^2} (G_2^{n+1} - g_2^{n+1}) d{\bf v}.
	\label{id5}
\end{equation}

Adding (\ref{id4}) and (\ref{id5}) and making use of the identity (\ref{id3}) we get
\begin{equation}
	3 \rho^{n+1} R T^{n+1} = \int_{\mathbb R^2} |{\bf v} - {\bf U}|^2 \tilde{g}_1^n d{\bf v} + \int_{\mathbb R^2} \tilde{g}_2^n d{\bf v}.
	\label{moment1b}
\end{equation}

Now, the parameters $\rho^{n+1}, {\bf U}^{n+1}$ and $T^{n+1}$ of $G_1^{n+1}$ and $G_2^{n+1}$ are given in terms of $\tilde{g}_1^n$ and $\tilde{g}_2^n$ from (\ref{moment1a}) and
(\ref{moment1b}). Hence the implicit steps (\ref{implicit_1}) and (\ref{implicit_2}) can be rewritten as
\begin{eqnarray}
	g_1^{n+1} = \frac{\tau \tilde{g}_1^n + \Delta t G_1^{n+1}}{\tau + \Delta t}\\
	g_2^{n+1} = \frac{\tau \tilde{g}_2^n + \Delta t G_2^{n+1}}{\tau + \Delta t}.
\end{eqnarray}

\subsection{Velocity discretization}
We consider $N_v$ velocity grid points and a uniform velocity grid of size $\Delta v = 2 v_{\rm max}/ (N_v-1)$
We assume that the distribution function is negligible for $|v| > v_{\rm max} $ and  discretize $[-v_{\rm max},v_{\rm max}]\times [-v_{\rm max},v_{\rm max}] \times [-v_{\rm max},v_{\rm max}]$.
That means, for each velocity direction, we have the discretization points
$v_j = -v_{\rm max} + (j-1)\Delta {v}, j = 1,\ldots , N_v$. 	
The velocity cutoff $v_{\rm max}=C_{\rm cut}\,\sigma$ must satisfy
$C_{\rm cut}\ge 5$: at $C_{\rm cut}=3.5$ the discrete Maxwellian quadrature
gives $T=269.3$\,K instead of the prescribed $270$\,K, and this mismatch
with the wall temperature seeds a slowly growing instability, while
$C_{\rm cut}=5.5$ is accurate to below $10^{-3}$\,K.

\subsection{Spatial discretization}

We discuss the spatial discretization and upwinding procedures for first-order schemes.

\subsubsection{Approximation of spatial derivatives}

In the above numerical schemes, an approximation of the spatial derivatives of $f$ or $g_1$ and $g_2$ is required. In this subsection, we describe a least-squares approximation of the derivatives on the moving point cloud based on so-called generalized finite differences; see \cite{SKT} and references therein.  A stabilizing procedure using upwinding discretization will be described in the following.

We consider a three-dimensional spatial domain $\Omega \cup \Gamma  \in \mathbb{R}^3$, where $\Gamma$ is the boundary. We first approximate
the boundary of the domain by a set of discrete points called boundary particles.   In the second step, we approximate the
interior of the computational domain using another set of interior points or interior particles. The sum of boundary and interior points gives the total number of points.
We note that the boundary conditions are applied on the boundary points. If the boundary moves, the boundary points also move together with the boundary. The initial generation of grid points
can be regular as well as arbitrary.  When the points
move, they can form a cluster or can scatter away from each other. In these cases, either some grid points have to be removed, or new grid points have to be added. We will describe this particle management in a forthcoming subsection.

Let $f(\bf {x})$ be a scalar function and $f_i$ its discrete values at ${\bf x}_i = (x_i, y_i, z_i)$. We consider the problem of approximating the spatial derivatives of a function $ f $ at ${\bf x}_i$ from the values of its neighboring points. In principle, all points can be neighbors; however,  we associate a weight function such that nearby particles have higher and
far away particles have lesser influence, such that we obtain an accurate approximation of derivatives. So, one can choose a weight as a function of the distance between the central point and its neighbor.  The distance function decays as the distance goes to infinity. In this paper, we have considered a Gaussian function, but other choices are possible (see, for example, \cite{sonar, Yom} for other classes of weight functions). In order to limit the number of neighboring points, we consider only the neighbors inside a circle of radius $h$ with center ${\bf x}_i$. We choose the radius $h$ to be some factor of the average spacing $\Delta x$,  such that we have at least a minimum number of neighbors for the least-squares approximation,  even next to the boundary. In the case of a regular grid far from the boundary, one might consider using smaller values of $h$. Such adaptive choice of  $h$ has been considered, for example, in \cite{kuhnert}. For the sake of simplicity, we have chosen a constant $h = 3.1~\Delta x$ in this paper, which gives a sufficiently large number of neighbours.
Let $P({\bf x}_{i_j}; h)=\{ {\bf x}_{i_j}, j=1,\ldots, m(h) \} $ denote the set of neighbor points of $ {\bf x}_i$ inside the disc of radius $h$. We note that the number $m$ of neighbours depends on  $h$.
In all calculations, we have considered the following truncated Gaussian weight function
\begin{eqnarray}
	w( \| {\bf x}_{i_j} - { \bf x}_i \|; h) =
	\left\{
	\begin{array}{l}
		\exp (- \alpha \frac{\|  {\bf x}_{i_j} - { \bf x}_i \|^2  }{h^2} ),
		\quad \mbox{if}  \frac{ \|  {\bf x}_{i_j} - { \bf x}_i \|} {h} \le 1
		\\ \nonumber
		0,  \qquad \qquad \mbox{else},
	\end{array}
	\right.
	\label{weight}
\end{eqnarray}
with $ \alpha $ a user-defined positive constant, chosen here as $\alpha = 6$, so that the influence of far neighbor grid points is negligible.
This choice is based on previous experiences \cite{kuhnert, TKH09}.

We consider the $m$ Taylor expansions of $f( x_{i_j}, y_{i_j}, z_{i_j})$ around $( x_i, y_i, z_i)$
\begin{equation}
	f_{i_j} = f( x_{i_j}, y_{i_j}, z_{i_j}) =  f_i  + \frac{\partial f}{\partial x} ( x_{i_j} - x_i) + \frac{\partial f}{\partial y} ( y_{i_j} -  y_i) + \frac{\partial  f}{\partial z} ( z_{i_j} -  z_i)  + e_{i_j},
	\label{taylor}
\end{equation}
for $j=1, \ldots, m$, where $e_{i_j}$ is the error in the Taylor's expansion. We note that the discrete values $ f_{i_j}$ are known. The system of equations (\ref{taylor}) can be re-expressed as
\begin{eqnarray}
	f_{i_1} - f_i &=&  \frac{\partial f}{\partial x}( x_{i_1} - x_i) + \frac{\partial f}{\partial y}( y_{i_1} -  y_i) + \frac{\partial f}{\partial z}( z_{i_1} -  z_i)  + e_{i_1} \nonumber \\
	\vdots & = & \vdots \\
	f_{i_m} - f_i &=&  \frac{\partial f}{\partial x}( x_{i_m} -  x_i) + \frac{\partial  f}{\partial y}( y_{i_m} -  y_i) +  \frac{\partial f}{\partial z}( z_{i_m} -  z_i) + e_{i_m}.  \nonumber
\end{eqnarray}
The system of equations can be written in  matrix form as
\begin{equation}
	{\bf e } = {\bf b} - M {\bf a},
\end{equation}
where ${\bf e} = [e_{i_1}, \ldots, e_{i_m}]^T$, ${\bf a} = [\frac{\partial f}{\partial x}, \frac{\partial f}{\partial y}, \frac{\partial f}{\partial z}]^T,
{\bf b} =[ f_{i_1} -  f_i \ldots,  f_{i_m} -  f_i]^T$
and
\begin{eqnarray}
	M=\left( \begin{array}{ccc}
		x_{i_1}-  x_i &  y_{i_1} -  y_i &  z_{i_1} -  z_i  \\
		\vdots & \vdots & \vdots \\
		x_{i_m} -  x_i &  y_{i_m} -  y_i & z_{i_m} -  z_i
	\end{array} \right)
	\nonumber.
\end{eqnarray}
If the number of neighbors $m$ is larger than three, this system of equations is overdetermined for three unknowns ${\bf a} = [\frac{\partial f}{\partial x}, \frac{\partial f}{\partial y}, \frac{\partial f}{\partial z}]^T$. The unknowns
$\bf{a}$ are obtained from the weighted least squares method by minimizing the quadratic form
\begin{equation}
	J = \sum_{j=1}^m w_j e_j^2 = (M{ \bf a} - {\bf b})^T W (M{\bf a } - {\bf b}),
\end{equation}
where $W=w_j\delta_{j k}, k = 1,\ldots, m $ is the diagonal matrix. The minimization of $J$ yields
\begin{equation}
	{\bf a} = (M^TWM)^{-1}(M^TW) ~{\bf b}.
	\label{solution}
\end{equation}
Denoting $S = (S_{ij}) = (M^TWM)^{-1}$ as the $3\times 3$ symmetric matrix and denoting $dx_j =  x_{i_j} -  x_i, dy_j =  y_{i_j} -  y_i, dz_j =  z_{i_j} -  z_i $ the vector $(M^TW)~{\bf b}$ is explicitly given by
\begin{eqnarray}
	(M^TW)~{\bf b} = \left( \begin{array}{c}
		\sum_{j=1}^m w_j dx_j ( f_{i_j} -  f_i)  \\
		\sum_{j=1}^m w_j dy_j ( f_{i_j} -  f_i )  \\
		\sum_{j=1}^m w_j dz_j ( f_{i_j} - f_i)  \\
	\end{array} \right).
	\label{mls_rhs}
\end{eqnarray}

Reading off the components of (\ref{solution}) we obtain the approximation of the spatial derivatives
\begin{eqnarray}
	\frac{\partial f}{\partial x} &=& \sum_{j=1}^m w_j ( S_{11}~dx_j  + S_{12}~dy_j + S_{13}~dz_j)( f_{i_j} -  f_i )
	\label{x_derivative}
	\\
	\frac{\partial f}{\partial y} &=& \sum_{j=1}^m w_j (S_{21}~dx_j  + S_{22}~dy_j + S_{23}~dz_j)( f_{i_j} -  f_i )
	\label{y_derivative}
	\\
	\frac{\partial  f}{\partial z} &=& \sum_{j=1}^m w_j ( S_{31}~dx_j  + S_{32}~dy_j + S_{33}~dz_j)( f_{i_j} - f_i).
	\label{z_derivative}
\end{eqnarray}

The corresponding computations are done for the two-dimensional case.
We note that higher-order approximations are obtained by using higher-order Taylor expansions in  (\ref{taylor}).
We refer to  \cite{TKH09} for details.

In the above least-squares approximation, a function is approximated at an arbitrary point from its neighboring points, and the distribution of these points can be arbitrary.
Such a straightforward least-squares approximation leads to a central difference scheme. In case of  discontinuities in the solution, this will lead to  numerical oscillations, and one has to
introduce additional numerical viscosity.
This can be done in the least squares framework by a suitable upwind reconstruction described in the next subsection.

\subsubsection{First order upwind scheme}

First, we describe the procedure for two-dimensional physical space. The naive method for computing upwind derivatives is to only use neighbor particles in the 'upwind'-half of the plane from the central particle. For example, we compute the x-partial derivatives of $f$  in the following way. If $v_1 - U_1 > 0$, we sort the neighbors with $x_{i_j} \le x_i$ and then compute the derivatives considering Taylor's expansion (\ref{taylor}) and then using the least squares method.
Similarly, if $v_1-U_1 < 0$ we use the set of neighbors with $x_{i_j} \ge x_i$.  For $y$-derivatives, we use the upper half and lower half of the neighbors according to the sign of $v_2 - U_2$. In a two-dimensional case, we have to use four stencils, and in a three-dimensional case, we have to use 6 stencils, which is time-consuming. Moreover, the positivity of the scheme is not guaranteed. Therefore, we have used a positive scheme suggested in \cite{Praveen}.

The scheme suggested in \cite{Praveen} makes use of a single central stencil and is, thus, computationally more efficient than the approach described above. Moreover, due to the additional diffusion that is added, the scheme is positive. In the following, we describe briefly the positive upwind scheme reported in \cite{Praveen}.

{\bf Two-dimensional physical space:} In the explicit step of \ref{explicit}, we have to compute the flux, for example, at the discrete point $x_i$
\begin{equation}
	{\bf Q}(f)_i =  (v_1 - U_1)\frac{\partial f}{\partial x} + (v_2 - U_2) \frac{\partial f}{\partial y}.
	\label{flux}
\end{equation}
Using (\ref{x_derivative} - \ref{y_derivative}) we can re-express (\ref{flux}) as
\begin{equation}
	{\bf Q}(f)_i = (v_1 - U_1) \sum_{j=1}^m \alpha_{ij}(f_{i_j} - f_i) + (v_2 - U_2)\sum_{j=1}^m \beta_{ij}(f_{i_j} - f_i),
\end{equation}
where
\begin{eqnarray*}
	\alpha_{ij}&=& w_j ( S_{11}~dx_j  + S_{12}~dy_j )
	\\
	\beta_{ij} &=& w_j (S_{21}~dx_j  + S_{22}~dy_j )
\end{eqnarray*}

Let $\phi_{ij} = \mbox{atan} 2\left({\frac{dy_j}{dx_j}}\right)$ be the angle between positive x-axis and the line ${x_ix}_{i_j}$ joining two points $x_i$ and $x_{i_j}$. Let ${\bf n}_{ij} = ( \cos\phi_{ij}, \sin\phi_{ij})$ be the unit vector along ${x_i x_{i_j}}$. Let ${\bf t}_{ij} = (-\sin\phi_{ij}, \cos\phi_{ij})$ be the unit vector orthogonal to ${\bf n}_{ij}$ so that ($({\bf n}_{ij}, {\bf t}_{ij})$ form a right handed coordinate system. Then the convective vector ${\bf v} - {\bf U}$ can be written in terms of this coordinate system after a rotational transformation
\begin{eqnarray*}
	v_1 - U_1 &=& ( ({\bf v} - {\bf U})\cdot {\bf n}_{ij}) \cos\phi_{ij} - (({\bf v} - {\bf U})\cdot {\bf t}_{ij})\sin\phi_{ij}\\
	v_2 - U_2 &=& ( ({\bf v} - {\bf U})\cdot {\bf n}_{ij}) \sin\phi_{ij} + (({\bf v} - {\bf U})\cdot {\bf t}_{ij})\cos\phi_{ij}.
\end{eqnarray*}

Then, after using mid-point flux approximations along ${\bf n}_{ij}$ and ${\bf t}_{ij}$ together with the upwind approximation of the flux, the upwind approximation of (\ref{flux}) is obtained in the form
\begin{eqnarray}
	{\bf Q}(f)_i = \sum_{j=1}^m \bar{\alpha}_{ij} \left[ ({\bf v} -{\bf U})\cdot {\bf n}_{ij} - |({\bf v} - {\bf U})\cdot {\bf n}_{ij}|\right] (f_{i_j} - f_i) + \nonumber
	\\
	\quad \quad \quad\quad \sum_{j=1}^m \bar{\beta}_{ij} \left[ ({\bf v} - {\bf U})\cdot {\bf t}_{ij} - \mbox{sign} (\bar{\beta}_{ij}) |({\bf v} - {\bf U})\cdot {\bf t}_{ij}|\right] (f_{i_j} - f_i),
\end{eqnarray}
where
\begin{eqnarray*}
	\bar{\alpha}_{ij} &=& \alpha_{ij} \cos\phi_{ij} + \beta_{ij} \sin\phi_{ij} \\
	\bar{\beta}_{ij} &= & - \alpha_{ij}\sin\phi_{ij} + \beta_{ij}\cos\phi_{ij}.
\end{eqnarray*}

{\bf Three-dimensional physical space:} As in the two-dimensional case, the velocity is re-expressed in terms of the component along $x_ix_{i_j}$ and the two normal directions. We define these vectors as
\begin{eqnarray*}
	A=\left( \begin{array}{ccc}
		n_{ij}(1) & n_{ij}(2) &n_{ij}(3) \\
		t_{ij}(1) & t_{ij}(2) &t_{ij}(3) \\
		b_{ij}(1) &b_{ij}(2) &b_{ij}(3) \\
	\end{array} \right)
	= \left( \begin{array}{ccc}
		\sin\theta_{ij}\cos\phi_{ij} & \sin\theta_{ij}\sin\phi_{ij} & \cos\theta_{ij} \\
		\cos\theta_{ij}\cos\phi_{ij}& \cos\theta_{ij}\sin\phi_{ij} & -\sin\theta_{ij} \\
		-\sin\phi_{ij}& ~ \cos\phi_{ij}& ~ 0
	\end{array} \right).
\end{eqnarray*}
where $\phi_{ij} = \mbox{atan2}\left(\frac{dy_j}{dx_j}\right), \; \theta_{ij} = \arccos\left(\frac{dz_j}{r}\right), \; r = \sqrt{dx_j^2 + dy_j^2 + dz_j^2}$.

With this convention, the rows of $A$ form an orthonormal triad whose
normal ${\bf n}_{ij}$ points along $(dx_j,dy_j,dz_j)/r$, so the upwind
selection $\left[\,({\bf v}-{\bf U})\cdot{\bf n}_{ij} -
|({\bf v}-{\bf U})\cdot{\bf n}_{ij}|\,\right]$ picks out exactly the
characteristics arriving from the direction of neighbour $j$.

Now using (\ref{x_derivative} - \ref{z_derivative}) we can re-express (\ref{flux}) as
\begin{equation}
	{\bf Q}(f)_i = (v_1 - U_1) \sum_{j=1}^m \alpha_{ij}(f_{i_j} - f_i) + (v_2 - U_2)\sum_{j=1}^m \beta_{ij}(f_{i_j} - f_i) 	+ (v_3 - U_3) \sum_{j=1}^m \gamma_{ij}(f_{i_j} - f_i),
\end{equation}
where
\begin{eqnarray*}
	\alpha_{ij}&=& w_j (S_{11}~dx_j  + S_{12}~dy_j + S_{13}~dz_j )
	\\
	\beta_{ij} &=& w_j (S_{21}~dx_j  + S_{22}~dy_j + S_{23}~dz_j )
	\\
	\gamma_{ij} &=& w_j(S_{31}~dx_j  + S_{32}~dy_j + S_{33}~dz_j )
\end{eqnarray*}

Writing the convective velocity vector ${\bf v} - {\bf U}$ in terms of coordinate systems $(n_{ij}, t_{ij}, b_{ij})$ with
\begin{eqnarray*}
	\left( \begin{array}{c}
		v_1 - U_1\\
		v_2 - U_2\\
		v_3 - U_3\\
	\end{array} \right) = A^T
	\left( \begin{array}{c}
		\langle {\bf v} - {\bf U}, {\bf n}_{ij} \rangle \\
		\langle {\bf v} - {\bf U} ,{\bf t}_{ij} \rangle\\
		\langle {\bf v} - {\bf U}, {\bf b}_{ij} \rangle\\
	\end{array} \right)
\end{eqnarray*}
and the rotational coefficients $\bar{\alpha}_{ij}, \bar{\beta}_{ij}, \bar{\gamma}_{ij}$ as
\begin{eqnarray*}
	\left( \begin{array}{c}
		\bar{\alpha}_{ij}\\
		\bar{\beta}_{ij}\\
		\bar{\gamma}_{ij}\\
	\end{array} \right) = A^t
	\left( \begin{array}{c}
		\alpha_{ij}\\
		\beta_{ij}\\
		\gamma_{ij}\\
	\end{array} \right)
\end{eqnarray*}
We have the upwind approximation of the flux, which is given as
\begin{eqnarray}
	{\bf Q}(f)_i = \sum_{j=1}^m \bar{\alpha}_{ij}  \left[ ({\bf v} - {\bf U})\cdot {\bf n}_{ij} - |({\bf v} - {\bf U})\cdot {\bf n}_{ij}|\right] (f_{i_j} - f_i) + \nonumber
	\\
	\quad \quad \quad\quad \sum_{j=1}^m \bar{\beta}_{ij} \left[ ({\bf v} - {\bf U})\cdot {\bf t}_{ij} - \mbox{sign}(\bar{\beta}_{ij}) |({\bf v} -{\bf U})\cdot {\bf t}_{ij}|\right] (f_{i_j} - f_i) + \nonumber
	\\
	\quad \quad \quad\quad \sum_{j=1}^m \bar{\gamma}_{ij} \left[ ({\bf v} -{\bf U})\cdot {\bf b}_{ij} - \mbox{sign}(\bar{\gamma}_{ij}) |({\bf v} - {\bf U})\cdot {\bf b}_{ij}|\right] (f_{i_j} - f_i)
	\label{upwind_3d}
\end{eqnarray}

In three dimensions, all velocity components are physical,
so the solver evolves the single distribution $f$ on the full
$(v_1,v_2,v_3)$ grid, and the projections in (\ref{upwind_3d}) use the
resolved node value $v_3$.  The transport update is
\begin{equation}\label{box3drel}
	\widetilde{f}^{\,n} \;=\; f^{n}
	\;-\; \Delta t\,{\bf Q}(f)_i^{n},
\end{equation}
with ${\bf Q}$ from (\ref{upwind_3d}), moments and
Maxwellian follow the three-dimensional relations of
Section~\ref{sec:model}, which are exact rather than
$O(\mathrm{Kn})$ approximations.

\subsection{Mass-conserving positivity fix}

The least squares reconstruction is not positivity
preserving. The transported distribution $\tilde g_1$ may undershoot to
negative values at isolated velocity nodes.  Before the moments are
evaluated, every particle applies a conservative clip-and-rescale.  Let
\begin{equation}
	m^{+} = \sum_{l:\,\tilde g_1^{(l)} \ge 0} \tilde g_1^{(l)}\,\Delta v^2,
	\qquad
	m^{-} = \sum_{l:\,\tilde g_1^{(l)} < 0} \tilde g_1^{(l)}\,\Delta v^2 ,
	\label{pos_masses}
\end{equation}
so that $m^{+}\!+m^{-}$ is the (physical) discrete mass of the
particle.  Provided $m^{+}\!+m^{-} > 0$, the corrected distribution is
\begin{equation}
	\hat g_k^{(l)} \;=\;
	\begin{cases}
		0, & \tilde g_1^{(l)} < 0,\\[2pt]
		s\, \tilde g_k^{(l)}, & \tilde g_1^{(l)} \ge 0,
	\end{cases}
	\qquad
	s \;=\; \frac{m^{+}+m^{-}}{m^{+}} \;\le\; 1,
	\qquad k = 1,2 .
	\label{pos_fix}
\end{equation}
Three properties follow directly.  (i) Positivity:
$\hat g_1 \ge 0$ at every node.  (ii) Exact mass conservation:
$\sum_l \hat g_1^{(l)}\Delta v^2 = s\,m^{+} = m^{+}\!+m^{-}$, the mass of
the uncorrected distribution.  (iii) Temperature consistency: both
layers are modified with the same node set and scale factor $s$,
so the ratio of the $g_2$- to the $g_1$-moment (with the
reduced-temperature contribution in (\ref{macroeq})) is preserved.
The fix is the identity whenever $\tilde g_1 \ge 0$, so accuracy is
untouched where the scheme is already positive, and the induced
perturbation of momentum and energy is bounded by the clipped mass
$|m^{-}|$, which vanishes in the same limit.  In three dimensions the
identical operator (\ref{pos_fix}) is applied to the single distribution
$f$, with the sums in (\ref{pos_masses}) running over the $N_v^3$ nodes
and $\Delta v^3$ in place of $\Delta v^2$.
\subsection{Particle management} In this section we discuss briefly some technical aspects of meshfree methods, such as, for example, neighbour searching, adding as well as removing particles; see \cite{kuhnert, DTKB08, TKR3} for more details.

\subsubsection{Neighbour search}
Searching neighbouring grid points at an arbitrary position is the most important and most   time-consuming
part of the meshfree method. Since the positions of grid points change over time, we have to search neighbor lists at every time step. We refer to \cite{TKR3} for details.

\subsubsection{Adding  and removing points}

When grid points move in time, they either scatter or cluster. When they scatter, there may be holes in the computational domain, and there are not enough neighbors at these points,  such that the numerical scheme becomes unstable. Therefore, one has to add new particles. On the other hand, if grid points cluster, there might be unnecessarily large numbers of neighbors at certain points, which requires more computational time and leads to a strongly non-uniform distribution of points, causing a bad condition number. In this case, if two points are close to each other, we remove both of them and introduce a new grid point at the midpoint. Then we update the distribution function $f(t,{\bf x}, {\bf v})$ on these new grid points with the help of the least squares method.
We refer to \cite{TKR} for a more detailed discussion.

\section{HPC Architecture and CUDA Implementation}
\label{sec:HPC}
High-performance computing (HPC) offers two complementary routes to parallelism, and this work uses both. In the distributed-memory model, a problem is divided among processes that each own part of the data in their own address space and cooperate through the explicit message passing interface (MPI). Because each process holds its own memory, the memory available to a calculation grows with the number of nodes, so problems too large for any single machine become reachable.  That capacity is only realised, however, if the data is genuinely partitioned. A parallelisation that merely divides the work while every process retains a
complete copy of the data leaves the aggregate memory of the machine unused,
and its cost is the communication needed wherever a computation reaches
beyond the data a process owns, which, for a meshfree scheme whose stencil
is defined by a geometric radius rather than by a fixed index pattern,
requires care at the partition boundaries. 

Modern GPUs, developed by NVIDIA and AMD, take the opposite approach: rather
than many independent address spaces, they provide thousands of smaller,
efficient cores sharing one high-bandwidth memory, which can run multiple
tasks simultaneously. Typical bottlenecks of GPUs are data transfer between
CPUs and GPUs, limited memory, and power consumption. Further, not all
computational tasks can be parallelized efficiently to utilize the GPU
architecture.

The two models therefore fail in
different ways: the distributed-memory machine by communication, the
accelerator by capacity and by transfer across the host interface. A
solver for the BGK equation meets both limits at once, since its
six-dimensional phase space makes memory rather than arithmetic the scarce
resource. The remainder of this section describes how each is addressed:
a slab domain decomposition with halo exchange for the distributed-memory
backend, and a node-parallel transport kernel with a fixed device footprint
for the GPU.

\subsection{HPC Architecture}

\textbf{Workstation:} The Windows implementation of the serial code is developed in Microsoft
Visual Studio 2026 (version 18) with the \texttt{/O2} optimization flag on a
dual-socket Intel(R) Xeon(R) Silver 4208 CPU @ 2.10\,GHz, providing 16
physical cores, with 192\,GB of RAM, under
Microsoft Windows Server 2019 Standard (version 1809). The distributed-memory implementation is built against Intel MPI, and the 16-rank configuration is measured on this same machine at one rank per physical core, so that the
parallel results and the sequential baseline share identical hardware.

\textbf{HPC:}  The
40-rank configuration is measured on a single compute node of Orchid with 40 logical cores (20 physical, two-way SMT), the institutional high-performance computing cluster at IIT Tirupati \cite{HPC}, carrying 128\, GB of RAM with jobs submitted through the PBS Pro scheduler.

Both distributed-memory configurations fit within a single node, so their
rank-to-rank exchanges are served by shared memory rather than by a network
interconnect. The communication cost reported below is therefore a lower
bound on what the same decomposition would incur across nodes, and these
results should be read as demonstrating that the slab decomposition removes
the per-node memory barrier rather than as a measurement of inter-node
efficiency.

\textbf{Google Colab:}		The GPU implementation is performed on a Google Colab Pro instance running Ubuntu 22.04.5 LTS (Jammy Jellyfish), equipped with an Intel(R) Xeon(R) CPU @ 2.20\, GHz with 12 virtual processors, 92\, GB RAM, and an NVIDIA A100-SXM4-80GB GPU (80\, GB HBM2 memory, 6912 CUDA cores, maximum GPU clock rate of 1.41\, GHz, and three copy engines). The CUDA implementation is developed using CUDA Toolkit v12.8.93 with the \texttt{-O3} and \texttt{-use\_fast\_math} compiler optimization flags.

\subsection{GPU Implementation: Data Layout and Host-Device Workflow}
\label{sec:cuda_layout}

The device cloud is stored in a structure-of-arrays layout: one contiguous
array per scalar field (positions, moments, wall data, neighbour lists) and
two large arrays \texttt{g} and \texttt{gt} holding each particle's
discrete distribution and its transport buffer.  With $N_v=20$ this is
$2N_v^2$ nodes per particle in 2D (the Chu pair, 12.8\,KB including the
buffer) and $N_v^3$ nodes in 3D (the single $f$, 64\, KB or 128 KB once the transport buffer is counted).  In the
benchmarks, $160{,}000$--$5{,}760{,}000$ particles are used in two
dimensions ($N_x=400$--$2400$) and $262{,}144$--$1{,}000{,}000$ in three
dimensions ($N_x=64, 80, 100$).

The host-device workflow is:
\begin{enumerate}
	\item The host allocates one device array per per-particle field, plus
	the two distribution arrays \texttt{g} and \texttt{gt}.
	\item The host generates the $N$ particles (indexed $0$ to $N-1$), the
	voxel decomposition of the domain, and, from the voxel neighbourhoods,
	the per-particle neighbour lists. These are rebuilt only if the cloud
	changes.
	\item The host initialises the distributions with the local Maxwellians
	and uploads the complete cloud to the device in a single pass.
	\item Time loop: each step executes the eight-kernel pipeline of
	Section~\ref{sec:cuda_upwind} entirely on the device. At every
	$50$ steps the moments $\rho,{\bf U},T$ are read back for
	output.
\end{enumerate}

\begin{figure}[!ht]
	\centering
	\begin{tikzpicture}[
		font=\small,
		node distance=4mm and 16mm,
		host/.style={rectangle, rounded corners=2pt, draw=black!75,
			fill=black!5, text width=40mm, align=center,
			minimum height=9mm},
		dev/.style={rectangle, draw=blue!55!black, fill=blue!7,
			text width=52mm, align=center, minimum height=7mm},
		seq/.style={-{Stealth[length=3mm]}, semithick, blue!55!black},
		xfer/.style={-{Stealth[length=3mm]}, thick, red!70!black},
		xlbl/.style={font=\scriptsize, red!70!black, align=center},
		note/.style={font=\scriptsize, black!60}
		]
		
		\node[dev]                  (mass1) {Global mass $m^n$\\[-1pt]
			{\scriptsize block reduction; $n_{\rm blocks}$ partials $\to$ host}};
		\node[dev, below=of mass1]  (copy)  {Transport-buffer prefill\\[-1pt]
			{\scriptsize \texttt{gt} $\leftarrow$ \texttt{g} (coalesced)}};
		\node[dev, below=of copy]   (trans) {\textbf{Semi-Lagrangian MLS transport}\\[-1pt]
			{\scriptsize one thread per (particle, velocity node):\\[-2pt]
				foot point $\to$ register-resident LS solve}};
		\node[dev, below=of trans]  (posi)  {Mass-conserving positivity fix};
		\node[dev, below=of posi]   (comm)  {Commit \texttt{g} $\leftarrow$ \texttt{gt}};
		\node[dev, below=of comm]   (mass2) {Global mass $m^{n+1}$ + rescale\\[-1pt]
			{\scriptsize enforces $m^{n+1}=m^{n}$ to round-off}};
		\node[dev, below=of mass2]  (wall)  {Diffuse-reflection wall BC};
		\node[dev, below=of wall]   (mom1)  {Moments $\rho,\ {\bf U},\ T$};
		\node[dev, below=of mom1]   (coll)  {Implicit BGK relaxation\\[-1pt]
			{\scriptsize one thread per (particle, velocity node)}};
		\node[dev, below=of coll]   (mom2)  {Moments (post-collision)};
		
		\foreach \a/\b in {mass1/copy, copy/trans, trans/posi, posi/comm,
			comm/mass2, mass2/wall, wall/mom1, mom1/coll,
			coll/mom2}
		\draw[seq] (\a) -- (\b);
		
		\draw[seq, dashed] (mom2.east) -- ++(7mm,0) |- (mass1.east)
		node[note, pos=0.25, rotate=90, yshift=2.5mm] {next step $\Delta t$};
		
		\node[host, left=of mass1] (init)
		{\textbf{CPU: initialise}\\[-1pt]
			{\scriptsize Generate particles and voxels,\\[-2pt]
				Neighbour lists, \\[-2pt] Initial Maxwellians}};
		\draw[xfer] (init.east) -- (mass1.west)
		node[xlbl, midway, above, sloped]
		{H$\,\to\,$D}
		node[xlbl, midway, below, sloped]
		{Once};
		
		\node[host, left=of mom2] (outp)
		{\textbf{CPU: output}\\[-1pt]
			{\scriptsize Diagnostics, VTK,\\[-2pt] Time series}};
		\draw[xfer, dashed] (mom2.west) -- (outp.east)
		node[xlbl, midway, above]
		{D$\,\to\,$H}
		node[xlbl, midway, below]
		{only $\rho,{\bf U},T$};
		
		\begin{scope}[on background layer]
			\node[draw=blue!45!black, dashed, rounded corners=3pt, inner sep=3mm,
			fit=(mass1)(mom2), label={[blue!45!black, font=\small]
				north: GPU - distribution resident on the device}] {};
		\end{scope}
	\end{tikzpicture}
	\caption{GPU execution model. The point cloud and the discrete
		distribution are uploaded once. The whole time step then runs as a
		sequence of kernels on device-resident data. The only recurring
		device-to-host traffic is the $n_{\rm blocks}$ partial sums of the mass
		reduction (kilobytes) and, at output steps, the macroscopic moments
		($\rho,{\bf U},T$, $5N$ doubles), the $N\!\cdot\!N_v^d$-sized distribution never
		leaves the GPU.}
	\label{fig:gpu_pipeline}
\end{figure}
\subsection{GPU Implementation: The Transport Step and Kernel Pipeline}
\label{sec:cuda_upwind}

The GPU code realises the explicit transport substep (\ref{explicit}) in an
equivalent semi-Lagrangian form.  Every velocity node traces its
characteristic foot point
${\bf x}_f = {\bf x}_i - ({\bf v}-{\bf U})\,\Delta t$ and the distribution
is reconstructed there with the weighted least squares apparatus of
Section \ref{sec:cuda_upwind}, the interpolation supplying the stabilising diffusion that the
upwind selection provides in the flux form.  One GPU thread owns one
velocity node of one particle ($N\cdot N_v^2$ threads in 2D,
$N\cdot N_v^3$ in 3D), so the fine-grained kinetic workload saturates the
device.  The per-step pipeline is:
\begin{enumerate}
	\item \textbf{Prefill:} the transport buffer is set to the current
	distribution, $\mathtt{gt}\leftarrow\mathtt{g}$, in one coalesced pass,
	so nodes the transport kernel skips (outgoing boundary nodes) retain
	their values.
	\item \textbf{Transport:} each thread writes the least squares
	reconstruction at its node's foot point into $\mathtt{gt}$.
	\item \textbf{Positivity fix:} negative nodes are clipped and the
	remainder is rescaled per particle by the mass-conserving correction
	(\ref{pos_fix}).
	\item \textbf{Commit:} $\mathtt{g}\leftarrow\mathtt{gt}$, completing the
	Jacobi update -- every node was reconstructed from the unmodified old
	distribution.
	\item \textbf{Global mass correction:} the total mass is summed by a
	two-stage block reduction and $\mathtt{g}$ is uniformly rescaled to the
	pre-transport mass, enforcing conservation to round-off.
	\item \textbf{Boundary condition:} at wall particles the incoming
	velocity half-space is replaced by the diffuse-reflection Maxwellian at
	the wall velocity and temperature.
	\item \textbf{Moments:} $\rho$, ${\bf U}$ and $T$ are evaluated per
	particle by quadrature over the velocity grid.
	\item \textbf{Implicit relaxation:} every node is updated pointwise by
	the L-stable formula
	$g^{n+1} = (\tau\,\tilde g + \Delta t\,G^{n+1})/(\tau+\Delta t)$,
	stable for any $\Delta t/\tau$.
\end{enumerate}
Figure~\ref{fig:gpu_pipeline} summarises the execution model: one upload, a
fully device-resident pipeline per step, and only kilobyte-scale readbacks
(mass partials, moments at output steps).

\textbf{Two-dimensional branch.}  Each thread transports both Chu layers of
its node $(v_1,v_2)$.  At boundary particles, outgoing nodes
($({\bf v}-{\bf U}_w)\cdot{\bf n}\ge 0$) keep their prefilled value and are
re-imposed by the boundary kernel, while incoming nodes use the linear
basis.  The weighted normal matrix (Gaussian weights, monomial basis
centred at the foot point -- quadratic with 6 terms in the interior, linear
with 3 at boundaries) is shared by both layers, so one in-register
elimination solves two right-hand sides at once.  The block reduction of
pipeline step 5 is essential rather than cosmetic: a single-address
\verb|atomicAdd| reduction serialises on the A100 and was measured at
$2274$\,ms against $0.68$\,ms for the block-parallel equivalent at
$N=160{,}000$.

\textbf{Three-dimensional branch.}  The single distribution $f$ is
transported on the full $N_v^3$ grid, the resolved node value $v_3$
entering the foot point directly.  Two choices dominate performance.
First, the kernel is templated on the polynomial size (10 monomials
interior, 4 at boundaries, as separate instantiations), keeping all least
squares arrays register-resident; a runtime-sized system spills to local
memory and was measured to slow the transport kernel by an order of
magnitude.  Second, the normal matrix is symmetric positive definite, so
only its upper triangle is assembled and a pivot-free Cholesky
factorisation solves it with no pivot search, row swaps, or divergent
branches inside the warp.

\subsection{Slab Decomposition and MPI Implementation}
The first distributed-memory backend replicates the point cloud: every rank
allocates the entire cloud and advances only its own contiguous index range,
after which a global \texttt{Allgatherv} per time step restores a consistent
copy everywhere. It is correct by construction and needed no change to the
kernels, since a rank holding every particle always finds the moving
least-squares stencil of any node it updates resident locally. Its
limitation is that only the work is distributed, never the storage. Per-rank
memory is independent of $P$, so the reachable problem size is node memory
divided by ranks per node, unaffected by how many nodes the machine has, and
falling as ranks are added since each demands another full copy. At
$N_v = 20$ in two dimensions the 16-rank workstation admits $N_x \approx 940$
and the 40-rank node, with two and a half times the cores but two thirds of
the memory, only $N_x \approx 490$. The communication scales no better, the
collective carrying the whole distribution rather than a boundary layer of
it. In practice, this backend reached $N_x = 1000$ and no further.

The slab decomposition \cite{Pekurovsky2012} partitions the lattice along its outermost generation
axis ($z$ in three dimensions, $y$ in two) as shown in
Figure~\ref{fig:slab}(a). That axis is chosen deliberately rather than for
geometric reasons: because particles are generated in that order, a slab is a
contiguous index range of the cloud as it already exists, so no reordering,
renumbering or scatter is needed, and the kernels, which already bound their
loops by the rank's chunk, run unmodified on the local array. The partition
alone is not sufficient, however. The circle drawn in Figure~\ref{fig:slab}(a)
is the moving least-squares support of a node in the topmost owned row of
rank~0, and it reaches across the cut into rows that rank~1 owns; updating
that node therefore requires data the rank does not hold. Each rank
consequently allocates its owned rows together with $n_{ghost}$ ghost rows on each
side, borrowed read-only from its two neighbours, giving the local array of
Figure~\ref{fig:slab}(b). Since the MLS support radius is a fixed multiple of
the node spacing, $h = 2.5\,\Delta x$, the halo is
$n_{ghost} = \lceil h/\Delta x\rceil + 1 = 4$ rows at every resolution, so its
relative cost falls under refinement instead of saturating. With the ghosts
in place, every stencil resolves locally again, including at the partition
boundaries, and the per-rank footprint becomes
\begin{equation}
	M_{\rm rank} \;=\; M_1\left(\frac{1}{P} + \frac{2n_{ghost}}{N_{\rm rows}}\right),
	\qquad
	M_{\rm aggregate} \;=\; M_1\left(1 + \frac{2n_{ghost}P}{N_{\rm rows}}\right),
	\label{eq:slabmem}
\end{equation}
where $M_1$ denotes the serial footprint of the whole cloud, $P$ the number
of ranks, and $N_{\rm rows}$ the number of lattice rows along the slab axis
($N_y$ in two dimensions, $N_z$ in three, and equal to $N_x$ for the cubic
grids used here).

The ghost rows go stale as soon as the neighbours update their own, so they
are refreshed once per time step by the two point-to-point exchanges of
Figure~\ref{fig:halo}. Only the ghost rows' distributions travel, and only to
the two immediate neighbours, so the volume per rank is independent of $P$
rather than proportional to it. The arrangement carries one structural
requirement: a neighbour must be able to serve a complete halo from rows it
owns, which forces the smallest owned slab to be at least $2n_{ghost} = 8$ rows tall
and therefore bounds the usable rank count by $P \le N_{\rm rows}/2n_{ghost}$. Where
that fails, the solver reports the condition and declines to decompose rather
than exchanging second-hand ghosts.

\begin{figure}[!ht]
	\centering
	
	\begin{tikzpicture}[
		font=\small,
		lbl/.style={font=\scriptsize},
		tin/.style={font=\tiny},
		own/.style={fill=blue!60!black},
		nbr/.style={fill=green!40!black},
		gho/.style={draw=orange!85!black, fill=orange!18, line width=0.5pt},
		ar/.style={-{Stealth[length=2mm]}, semithick, black!60}
		]
		\def\s{0.34}
		\def\d{0.072}
		
		\begin{scope}
			\foreach \j in {0,...,15}{\foreach \i in {0,...,8}{
					\ifnum\j<8 \fill[own] (\i*\s,\j*\s) circle (\d);
					\else \fill[nbr] (\i*\s,\j*\s) circle (\d);\fi}}
			\draw[dashed, thick, black!75] (-0.35,7.5*\s) -- (8*\s+0.35,7.5*\s);
			\draw[red!75!black, thick] (4*\s,7*\s) circle (2.5*\s);
			\fill[red!75!black] (4*\s,7*\s) circle (0.085);
			\node[tin, black!70, anchor=west] at (8*\s+0.42,3.5*\s) {rank 0};
			\node[tin, black!70, anchor=west] at (8*\s+0.42,11.5*\s) {rank 1};
			\node[lbl, red!75!black, anchor=west] at (8*\s+0.42,7.5*\s) {$r=2.5\,\Delta x$};
			\node[lbl, rotate=90] at (-0.78,7.5*\s) {slab axis ($z$ in 3-D, $y$ in 2-D)};
			\node[lbl] at (8*\s/2,-0.55) {(a) global lattice, cut into two slabs};
		\end{scope}
		
		\begin{scope}[shift={(6.2,0)}]
			\foreach \j in {0,...,7}{\foreach \i in {0,...,8}{
					\fill[own] (\i*\s,\j*\s) circle (\d);}}
			\foreach \j in {8,...,11}{\foreach \i in {0,...,8}{
					\draw[gho] (\i*\s,\j*\s) circle (\d);}}
			\draw[dashed, thick, black!75] (-0.35,7.5*\s) -- (8*\s+0.35,7.5*\s);
			\draw[red!75!black, thick] (4*\s,7*\s) circle (2.5*\s);
			\fill[red!75!black] (4*\s,7*\s) circle (0.085);
			\draw[thick, black!60, |-|]        (8*\s+0.3,-0.4*\s) -- (8*\s+0.3,7.4*\s);
			\draw[thick, orange!85!black, |-|] (8*\s+0.3,7.6*\s)  -- (8*\s+0.3,11.4*\s);
			\node[tin, black!70, anchor=west]        at (8*\s+0.42,3.5*\s) {owned};
			\node[tin, orange!85!black, anchor=west] at (8*\s+0.42,9.5*\s) {$h=4$};
			\node[tin, orange!85!black, anchor=west] at (8*\s+0.42,8.9*\s) {ghost};
			\node[lbl] at (8*\s/2,-0.55) {(b) the array rank 0 allocates};
		\end{scope}
	\end{tikzpicture}
	\caption{Slab decomposition. (a) The lattice is cut along its outermost
		generation axis into contiguous ranges, one per rank: filled circles are
		owned rows, blue for rank~0 and green for rank~1. The red circle is the MLS
		support, $r = 2.5\,\Delta x$, of the node at its centre. (b) Rank~0 therefore
		also stores $h = 4$ ghost rows (open orange) copied read-only from rank~1,
		after which every stencil it must evaluate closes within the local array.
		Eight owned rows are drawn, the thinnest slab the scheme admits, so the halo
		appears proportionally thicker here than at production resolutions.}
	\label{fig:slab}
\end{figure}
\vspace{14mm}
\begin{figure}
	\centering
	
	\begin{tikzpicture}[
		font=\small,
		lbl/.style={font=\scriptsize},
		tin/.style={font=\tiny},
		own/.style={fill=blue!60!black},
		nbr/.style={fill=green!40!black},
		gho/.style={draw=orange!85!black, fill=orange!18, line width=0.5pt},
		up/.style={-{Stealth[length=2.4mm]}, thick, blue!55!black},
		dn/.style={-{Stealth[length=2.4mm]}, thick, green!35!black}
		]
		\def\s{0.34}
		\def\d{0.072}
		\def\G{3.4}
		
		\foreach \j in {0,...,7}{\foreach \i in {0,...,5}{
				\fill[own] (\i*\s,\j*\s) circle (\d);}}
		\foreach \j in {8,...,11}{\foreach \i in {0,...,5}{
				\draw[gho] (\i*\s,\j*\s) circle (\d);}}
		\node[tin, black!70] at (2.5*\s,-0.42) {rank 0};
		
		\begin{scope}[shift={(\G,0)}]
			\foreach \j in {4,...,7}{\foreach \i in {0,...,5}{
					\draw[gho] (\i*\s,\j*\s) circle (\d);}}
			\foreach \j in {8,...,15}{\foreach \i in {0,...,5}{
					\fill[nbr] (\i*\s,\j*\s) circle (\d);}}
			\node[tin, black!70] at (2.5*\s,-0.42) {rank 1};
		\end{scope}
		
		\draw[up] (5*\s+0.14,5.5*\s) -- (\G-0.14,5.5*\s);
		\node[tin, blue!55!black, anchor=south] at (\G/2+0.9*\s,5.75*\s)
		{$h$ rows};
		\draw[dn] (\G-0.14,9.5*\s) -- (5*\s+0.14,9.5*\s);
		\node[tin, green!35!black, anchor=south] at (\G/2+0.9*\s,9.75*\s)
		{$h$ rows};
		
		\node[lbl, anchor=west] at (\G+6*\s+0.3,9.5*\s) {phase 1: downward};
		\node[lbl, anchor=west] at (\G+6*\s+0.3,5.5*\s) {phase 2: upward};
		\node[lbl] at (\G/2+2.5*\s,-0.95)
		{each rank sends $h$ owned rows, receives $h$ ghost rows};
	\end{tikzpicture}
	\caption{Halo exchange, once per time step. Two \texttt{Sendrecv} phases
		refresh the stale ghost rows (open orange) from the neighbours' owned rows
		(filled): rank~1 sends downward, rank~0 upward. Only the ghost rows travel,
		and only between immediate neighbours, so the volume per rank is independent
		of the rank count. Wall-bounded axes have no partner on the outer side.}
	\label{fig:halo}
\end{figure}

\subsection{HPC Memory Management}
\label{sec:gpu_memory}

Device memory is dominated by the two distribution arrays \texttt{g} and
\texttt{gt}, together
$2\times N\times(\mathrm{nodes\ per\ particle})\times 8$ bytes, the $O(N)$
per-particle field arrays are negligible in comparison.  All particle
management runs on the host, and the cloud is re-uploaded only when the
particle count changes, so no device-side scratch copies are needed -- in
contrast to a device-resident particle manager, which requires up to three
full copies of the particle array.  Memory requirements are summarised in
Tables~\ref{table2d_mem} and~\ref{table3d_mem}.  Before the time loop, the
solver prints the GPU architecture and a per-allocation memory budget with
the predicted free memory, so an out-of-memory configuration is visible
before any kernel runs, and every allocation and launch is wrapped in a
checking macro that aborts with the failing file and line instead of
continuing with a null pointer (which would manifest as $\rho=0$ in all
subsequent moments).

The practical consequence of (\ref{eq:slabmem}) is set out in
Table~\ref{tab:slablimits}. In two dimensions, the halo is a negligible
fraction of a slab: at $N_v = 20$ the 16-rank workstation reaches
$N_x = 3860$, and the 40-rank HPC reaches $N_x=2958$. In three dimensions, a lattice of $N_{\rm rows}$ rows admits at most
$N_{\rm rows}/2n_{ghost}$ slabs, so neither the workstation nor the HPC can be used at full width.
Within the admissible range, the cost is no longer slight: on the workstation
the eight-rank
limit is $N_x = 125$ (Table \ref{tab:slablimits}), so the halo costs thirteen per cent for eightfold parallelism, whereas in two dimensions it costs under two. Four ghost rows are negligible against a two-dimensional slab
and ruinous against a three-dimensional one, so resolution and parallelism
now trade directly against one another.

\begin{table}[H]
	\begin{center}
		\begin{tabular}{|r|r|c|}
			\hline
			$N_x$ & $N=N_x^2$\rule[-2.0ex]{0pt}{4.6ex} & Device memory  \\ \hline
			$1000$ & 1,000,000  & 12.47\,GB \\
			$1200$ & 1,440,000 & 17.97\,GB \\
			$1600$ & 2,560,000 & 31.97\,GB \\
			$2000$ & 4,000,000 & 49.96\,GB \\
			$2400$ & 5,760,000 & 71.96\,GB \\                   
			\hline
		\end{tabular}
		\caption{GPU memory for two-dimensional problems ($N_v=20$,
			Chu pair $(g_1,g_2)$ on $2N_v^2$ nodes per particle): twice the
			distribution storage ($\mathtt{g}$ + $\mathtt{gt}$) plus the small
			per-particle field arrays.  Up to $1160^2$ fits a 16\,GB card whereas the
			$1600^2$ case (31.97\,GB) requires the 40\,GB A100. $2000^2$ case (49.96\,GB) requires an 80\, GB A100 card. $2500^2$ requires H200 or B200 cards.}
		\label{table2d_mem}
	\end{center}
\end{table}

\begin{table}[H]
	\begin{center}
		\begin{tabular}{|r|r|c|c|}
			\hline
			$N_x$ & $N=N_x^3$\rule[-2.0ex]{0pt}{4.6ex} & Device memory ($N_v=20$) & Device memory ($N_v=16$)\\ \hline
			$80$ & 512,000 & 61.33\,GB & 31.55\,GB\\
			$100$ & 1,000,000 & 119.21\,GB & 61.61\, GB \\
			$120$ & 1,728,000 & 205.99\,GB & 106.47\, GB \\
			\hline
		\end{tabular}
		\caption{GPU memory for three-dimensional problems
			(single distribution $f$ on $N_v^3$ nodes). The $80^3$ with $N_v=20$ case approaches the capacity of an
			80\,GB A100 GPU  whereas it can fit in 40\, GB A100 GPU for $N_v=16$. The $100^3$ with $N_v=16$ case easily fits in an 80\, GB A100 card, whereas $120^3$ requires H200 or B100 cards.}
		\label{table3d_mem}
	\end{center}
\end{table}

\begin{table}[H]
	\begin{center}
		\begin{tabular}{|r|r|c|c|c|c|}
			\hline
			$N_x$ &$N=N_x^2$\rule[2.0ex]{0pt}{1.6ex} & \multicolumn{2}{c|}{Workstation (16 ranks)} & \multicolumn{2}{c|}{HPC (40 ranks)}\\
			\cline{3-4}\cline{5-6}
			&  & $M_{rank}$ & $M_{aggregate}\rule[-1.0ex]{0pt}{2.6ex}$ & $M_{rank}$ & $M_{aggregate}$\\ \hline
			$1000$ &  1,000,000 &  0.84\,GB &  13.45\,GB & 0.39\,GB &  15.74\,GB\\
			$2000$ &  4,000,000 &  3.17\,GB &  50.74\,GB & 1.38\,GB &  55.31\,GB\\
			$3000$ &  9,000,000 &  6.99\,GB & 111.87\,GB & 2.97\,GB & 118.73\,GB$^\dagger$\\
			$4000$ & 16,000,000 & 12.30\,GB & 196.84\,GB$^\dagger$ & 5.15\,GB & 205.99\,GB$^\dagger$\\
			\hline
		\end{tabular}
		\caption{Distributed-memory footprint for the two-dimensional problems ($N_v=20$)
			under the slab decomposition. $^\dagger$ exceeds the memory available to that configuration ($M_1=183.52$\, GB on the 16-rank workstation, $M_1=115.6$\, GB on the 40-rank node).}
		\label{table2d_mpi}
	\end{center}
\end{table}

\begin{table}[H]
	\begin{center}
		\begin{tabular}{|r|r|c|c|c|c|}
			\hline
			$N_x$ & $N=N_x^3$\rule[2.0ex]{0pt}{1.6ex} & \multicolumn{2}{c|}{$N_v=20$} & \multicolumn{2}{c|}{$N_v=16$}\\
			\cline{3-4}\cline{5-6}
			& & $M_{rank}$ & $M_{aggregate}\rule[-1.0ex]{0pt}{2.6ex}$ & $M_{rank}$ & $M_{aggregate}$\\ \hline
			$80$  &   512,000 & 13.73\,GB & 109.9\,GB &  7.03\,GB &  56.3\,GB\\
			$100$ & 1,000,000 & 24.44\,GB & 195.5\,GB & 12.51\,GB & 100.1\,GB\\
			$120$ & 1,728,000 & 39.48\,GB & 315.9\,GB & 20.21\,GB & 161.7\,GB\\
			\hline
		\end{tabular}
		\caption{Distributed-memory footprint for three-dimensional problems under the slab
			decomposition at $P=8$. At $N_v=16$ the largest cases reachable are
			$125^3$ on the workstationand $105^3$
			on the HPC.}
		\label{table3d_mpi}
	\end{center}
\end{table}

\begin{table}[H]
	\centering
	\caption{Largest resolution reachable on a single A100 GPU and on the two
		distributed-memory configurations. The GPU figures are limited by device
		memory on one card, with no decomposition. For the CPU configurations, the
		full rank count is used in two dimensions, where the parallel limit lies
		within a few per cent of the serial one; in three dimensions the bound
		$P \le N_{\rm rows}/2n_{ghost}$ caps the usable rank count at 13 and 12
		respectively, so the decomposed figures are quoted for $P=8$ ranks.}
	\label{tab:slablimits}
	\begin{tabular}{|l|c|c|c|}
		\hline
		Configuration & $N_v$ & $N_x$ (2D), full rank & $N_x$ (3D), $P=8$\\
		\hline
		A100 80\,GB (single GPU) & 16 & 3072 & 105\\
		& 20 & 2457 &  84\\
		\hline
		16 ranks (workstation) & 16 & 4840 & 125\\
		& 20 & 3860 &  97\\
		\hline
		40 ranks (HPC)         & 16 & 3735 & 105\\
		& 20 & 2958 &  81\\
		\hline
	\end{tabular}
\end{table}

\subsection{CFL Condition and Adaptive Time Step}
\label{sec:timestep}
The semi-Lagrangian transport has no advective CFL
stability limit. The restriction is that each node's foot point
${\bf x}_f = {\bf x}_i - ({\bf v}-{\bf U})\Delta t$ stays inside its
particle's least-squares stencil, so the reconstruction interpolates
rather than extrapolates as $$\Delta t \le C_{\rm CFL}\,\frac{\Delta x_{\min}}{v_{\max}}, ~~~~C_{\rm CFL}=1.5$$  Violations are clamped to the safer
$C_{CFL}=0.75$, since a foot point at the stencil edge has a one-sided,
near-singular reconstruction (in the 3-D cavity
with a rigid body at $N_x=80$, the reference step gives $C_{\rm CFL}=1.54$, so it is reduced to the $C_{\rm CFL} = 0.75$ value, $7.29\times10^{-12}$\,s).

The step is selected adaptively, but only at start-up rather than at every step. $v_{\max}$ in the bound is the fixed extent of the velocity grid
which is the fastest characteristic that can exist (not the instantaneous
flow speed), so the admissible $\Delta t$ is invariant in time and is
evaluated (and clamped if necessary) only at start-up.  Per-step
re-evaluation would be needed only with an adaptive velocity grid.

The transport kernel dominates the arithmetic cost. Each phase-space node
needs $m$ weighted neighbour contributions to the normal equations and one
$p \times p$ least-squares solve. With $m \approx 20$, $p = 6$ in two
dimensions and $m \approx 65$, $p = 10$ in three, the floating-point
operations per time step, $F_{\rm 2D}$ and $F_{\rm 3D}$, are
\begin{equation}
	F_{\rm 2D} \approx 1.2\times10^{3}\, N_x^{2} N_v^{2},
	\qquad
	F_{\rm 3D} \approx 4.7\times10^{3}\, N_x^{3} N_v^{3},
	\label{eq:flops}
\end{equation}
scaling with the same $N_x^{d} N_v^{d}$ product as the memory footprint.
Since $\Delta t \propto \Delta x = L/(N_x-1)$, the step count at fixed
$t_{\rm final}$ grows as $N_{\rm steps} \propto N_x$, so the total work of a
complete run, $W = N_{\rm steps}F$, scales as
\begin{equation}
	W_{\rm 2D} \propto N_x^{3} N_v^{2},
	\qquad
	W_{\rm 3D} \propto N_x^{4} N_v^{3},
	\label{eq:work}
\end{equation}
and doubling $N_x$ costs a factor of eight in two dimensions and sixteen in
three. The A100 sustains $680$\,GFLOP/s in two dimensions and $2970$ in three
($7.0$ and $30.6\%$ of double-precision peak), the difference reflecting that
the neighbour list is reused across $N_v^{3}$ velocity nodes rather than
$N_v^{2}$. Indirect neighbour access, not arithmetic, is the limit, which is
what makes the layout of Section~\ref{sec:cuda_layout} load-bearing.

Table~\ref{table_dt} lists, for each resolution, the largest admissible time
step $\Delta t_{\max} = C\,\Delta x / v_{\max}$ with $\Delta x = L/(N_x-1)$
and the corresponding number of steps. The value $C = 1.5$ is the hard
stencil bound, at which the characteristic foot point lies on the edge of the
least-squares support, and the reconstruction becomes one-sided; all
computations reported here therefore use $C \le 0.75$. Since
$\Delta t_{\max} \propto \Delta x \propto (N_x-1)^{-1}$ while the cost per
step grows as $N_x^{d}$ at fixed $N_v$, the total work at fixed
$t_{\rm final}$ scales as $N_x^{3}$ in two dimensions and $N_x^{4}$ in three.
The uniform step $\Delta t = 10^{-13}$\,s used for the backend comparison of
the next section is far more conservative still: at $N_x = 100$ it
corresponds to $v_{\max}\Delta t/\Delta x = 0.013$, some $58$ times below the
recommended step. It is retained only so that every backend and problem is
timed on an identical step count. Production runs, which carry no such
constraint, use the $C \le 0.75$ step and are correspondingly cheaper.

\begin{table}
	\begin{center}
		\begin{tabular}{|r|r|c|c|r|c|r|}
			\hline
			\multirow{2}{*}{$N_x$} & \multirow{2}{*}{$N$} & \multirow{2}{*}{$\Delta x$ (m)}
			& \multicolumn{2}{c|}{Stencil bound (CFL $1.5$)}
			& \multicolumn{2}{c|}{Recommended (CFL $0.75$)} \\ \cline{4-7}
			& & & $\Delta t_{\max}$ (s) & steps & $\Delta t$ (s) & steps \\ \hline
			\multicolumn{7}{|c|}{Two dimensional} \\ \hline
			400  & 160\,000    & 2.506e-09 & 2.885e-12 & 1733 & 1.443e-12 & 3466 \\
			800  & 640\,000    & 1.252e-09 & 1.441e-12 & 3470 & 7.204e-13 & 6941 \\
			1200 & 1\,440\,000 & 8.340e-10 & 9.601e-13 & 5208 & 4.801e-13 & 10415 \\
			1600 & 2\,560\,000 & 6.254e-10 & 7.199e-13 & 6945 & 3.600e-13 & 13890 \\ 
			2000 & 4\,000\,000 & 5.003e-10 & 5.759e-13 &  8682 & 2.879e-13 & 17365 \\
			2400 & 5\,760\,000 & 4.168e-10 & 4.799e-13 & 10420 & 2.399e-13 & 20839 \\ \hline
			\multicolumn{7}{|c|}{Three dimensional} \\ \hline
			80 & 512\,000  & 1.266e-08 & 1.457e-11 & 343 & 7.286e-12 & 686 \\
			100 & 1\,000\,000 & 1.010e-08 & 1.163e-11 & 430 & 5.814e-12 &  860 \\
			120 & 1\,728\,000 & 8.403e-09 & 9.674e-12 & 517 & 4.837e-12 & 1034 \\          
			\hline
		\end{tabular}
		\caption{Foot-point time-step limit and step counts
			($t_{\rm final} = 5\times10^{-9}$\,s).}
		\label{table_dt}
	\end{center}
\end{table}

\section{Numerical results}
\label{sec:numerics}

As numerical test cases, we consider the kinetic-driven cavity with and without a rigid body, in two- and three-dimensional settings.

\subsection{Two-dimensional driven cavity problem}

In this case, we use the two-dimensional BGK model that makes use of the Chu-reduction (see section \ref{1st_order_step_red}).  We consider a square cavity  $[0, L] \times [0, L] $ with  $L=1\times10^{-6} m$. Here $L$ is considered as the characteristic length. We consider  Argon gas with diameter $d = 0.368\times 10^{-9} m$ , the Boltzmann constant $k_B = 1.3806\times 10^{-23} JK^{-1}$ and the gas constant $R=208 Jkg^{-1}K^{-1}$.
Initially, the gas is in thermal equilibrium, and the distribution function is the Maxwellian distribution with initial parameters ${\bf U}_0 = 0 ms^{-1}, T_0 = 270 K$. We consider densities  $\rho_0 = 1$ and $0.1$ such that their corresponding Knudsen numbers are $Kn = 1.1$ and $Kn = 0.11$, respectively. From Eq. \ref{tau} the corresponding relaxation times are $\tau = 3.7142\times10^{-10}$ and $\tau = 3.7142\times 10^{-9}$. We apply diffuse reflection boundary conditions on all walls, where we keep the constant temperature $T_{wall} = T_0$ and the wall velocity ${\bf U}_{wall}$. The top wall has the velocity ${\bf U}_{wall} = (1,0) ms^{-1}$ and the other walls have zero velocity. The number of velocity grid points is $N_v = 20$.
The spatial grid points are initially generated regularly with numbers ranging from 400 to 2400 for each spatial direction.
When they move, their distribution will become random after a few time steps. 
\begin{figure}[!ht]
	\centering
	
	\includegraphics[keepaspectratio=true, angle=0, width=0.48\textwidth]{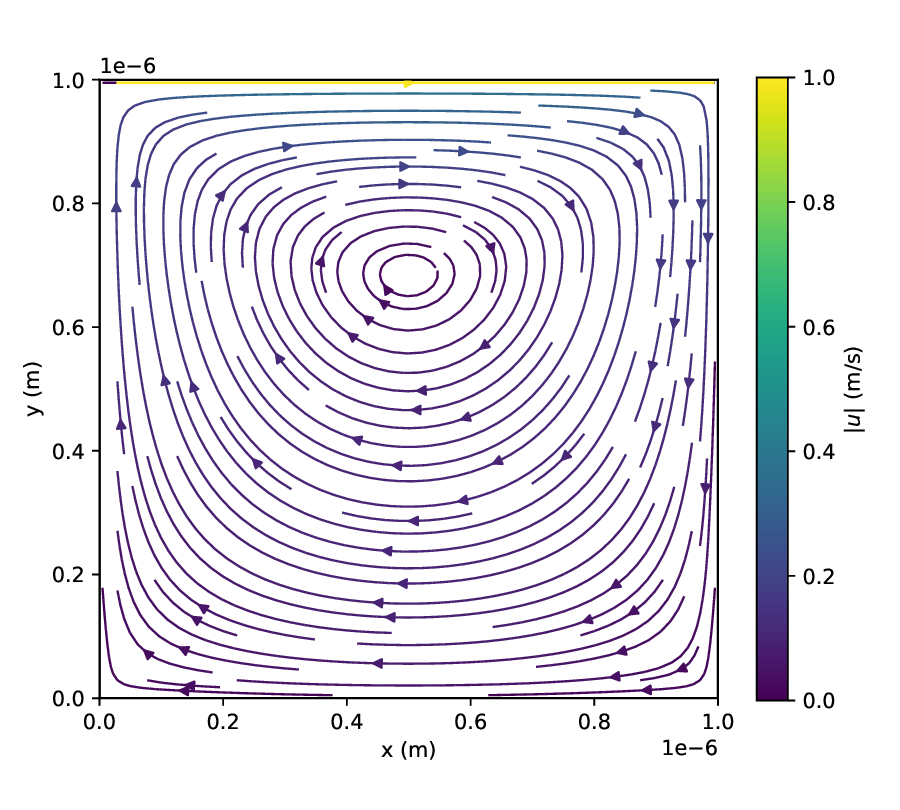}
	\includegraphics[keepaspectratio=true, angle=0, width=0.48\textwidth]{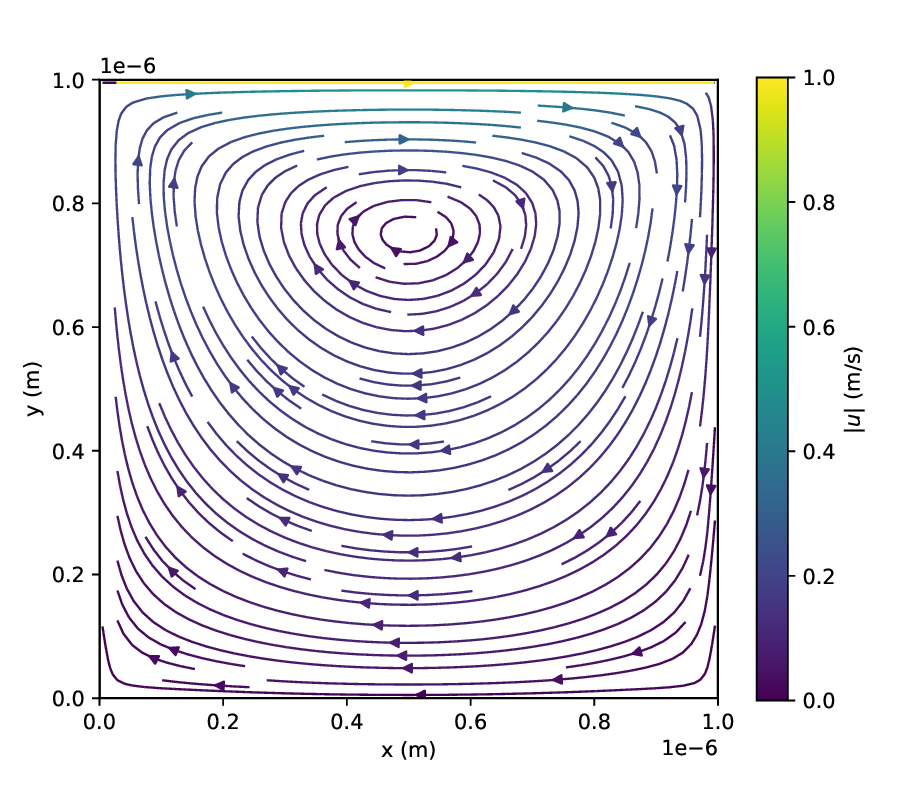}
	\caption{Two-dimensional driven cavity velocity vectors at $Kn = 1.1, Re = 2.968 \times 10^{-4}$ (left) and $Kn = 0.11, Re = 2.968 \times 10^{-5}$ (right) at $t=5\times 10^{-9}$.}
	\label{2dcavity}
	\centering
\end{figure}

The GPU implementation of the two-dimensional transport kernel is validated by monitoring the domain-averaged density $\bar\rho$ and temperature $\bar T$ over time.  With velocity cutoff $C_{\rm cut}=5.5\,\sigma$, the deviation of $\bar\rho$ from unity remains below $0.01\%$ and that of $\bar T$ from $T_0=270$\, K below $0.1$\, K throughout the simulation, confirming that the mass-conserving positivity fix and the Gaussian weight normalisation are correct.  


In Fig. \ref{2dcavity} we have plotted the velocity vectors in steady state for both Knudsen numbers. We note that the Mach and  Reynolds numbers are very small for these computations.
We have performed the same simulation on both the CPU and the GPU. In Table \ref{tab:bench2dcav} we present the GPU computation time and speedup compared to serial and parallel CPU computation times. As expected, we observe that the GPU speedup is higher compared to both serial and parallel CPUs when finer resolutions are used. For four million particles, the GPU is roughly $757$ times faster than a
single CPU core, $51$ times faster than $16$ MPI processes on the
workstation, and $29$ times faster than $40$ processes on the HPC node. The GPU completes the $2400 \times 2400$ case in approximately $56$ hours,
against $117$ days on the workstation and $69$ days on the HPC node. With the adaptive time step, the GPU takes about a day, while the workstation and the HPC node take 49 and 29 days, respectively.

\begin{table}[H]
	\begin{center}
		\begin{tabular}{|r|r|c|c|c|c|c|c|}
			\hline
			\multirow{2}{*}{$N_x$} & \multirow{2}{*}{$N$} & \multicolumn{4}{c|}{Time per step (s)} & GPU full run \\ \cline{3-6}
			& & Serial & Workstation & HPC & GPU  & (h) \\ \hline
			400  & 160,000    &  85.54 & 6.32 & 4.4 & 0.12 & 1.78 \\
			800  & 640,000    & 343.96  & 23.56 & 14.2 & 0.45 & 6.3 \\
			1200 & 1,440,000 & 762.55 & 52.59   & 30.14 & 1.03 & 14.2 \\
			1600 & 2,560,000 & 1395.05  & 92.41   & 54.34& 1.81 & 25.1 \\   		
			2000 & 4,000,000 & 2085.76  & $142.31$   & 82.87 & 2.81 & 39 \\
			2400 & 5,760,000 & 3049.84 & $202.63$   & 119.24 & 4.03 & 56 \\
			\hline
		\end{tabular}
		\caption{2D driven cavity: measured time per step for the serial, MPI
			and A100 GPU on Google Colab Pro ($\Delta t=10^{-13}, t_{\rm final}=5\times 10^{-9}$, 50,000 steps).}
		\label{tab:bench2dcav}
	\end{center}
\end{table}

\subsection{Three-dimensional driven cavity problem}

\begin{figure}
	\centering
	\includegraphics[keepaspectratio=true, angle=0, width=0.48\textwidth]{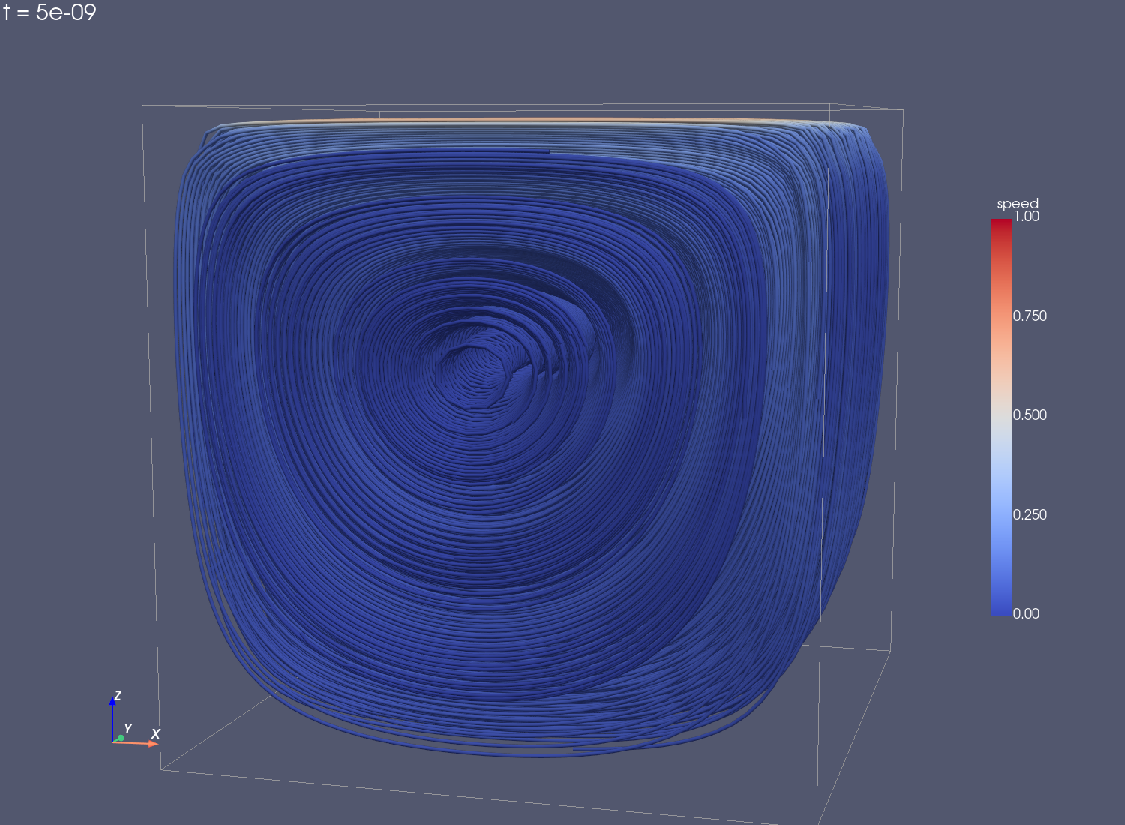}
	\includegraphics[keepaspectratio=true, angle=0, width=0.48\textwidth]{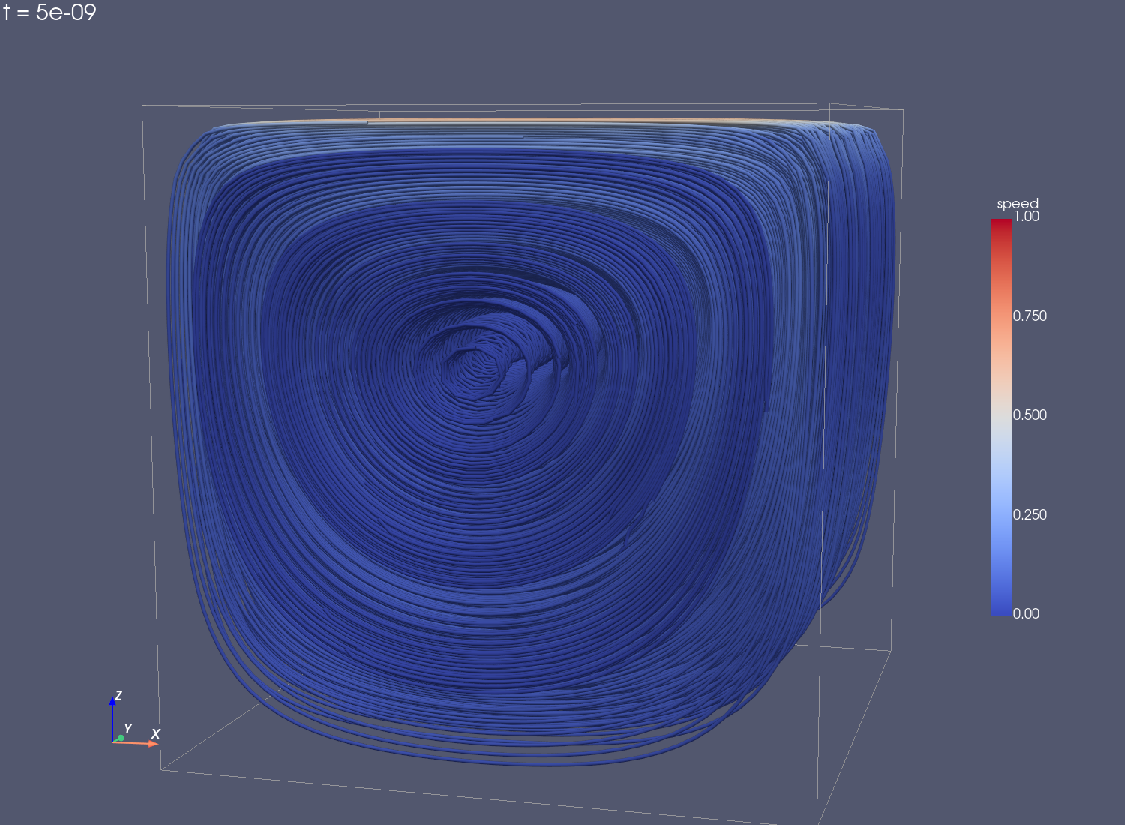}
	\caption{Three-dimensional driven cavity stream plot for Knudsen numbers and Reynolds numbers $Kn = 1.1, Re = 2.968 \times 10^{-4}$ (left) and $Kn = 0.11, Re = 2.968 \times 10^{-5}$ (right). }
	\label{3dcavity}
	\centering
\end{figure}

\begin{figure}
	\centering
	\includegraphics[keepaspectratio=true, angle=0, width=0.48\textwidth]{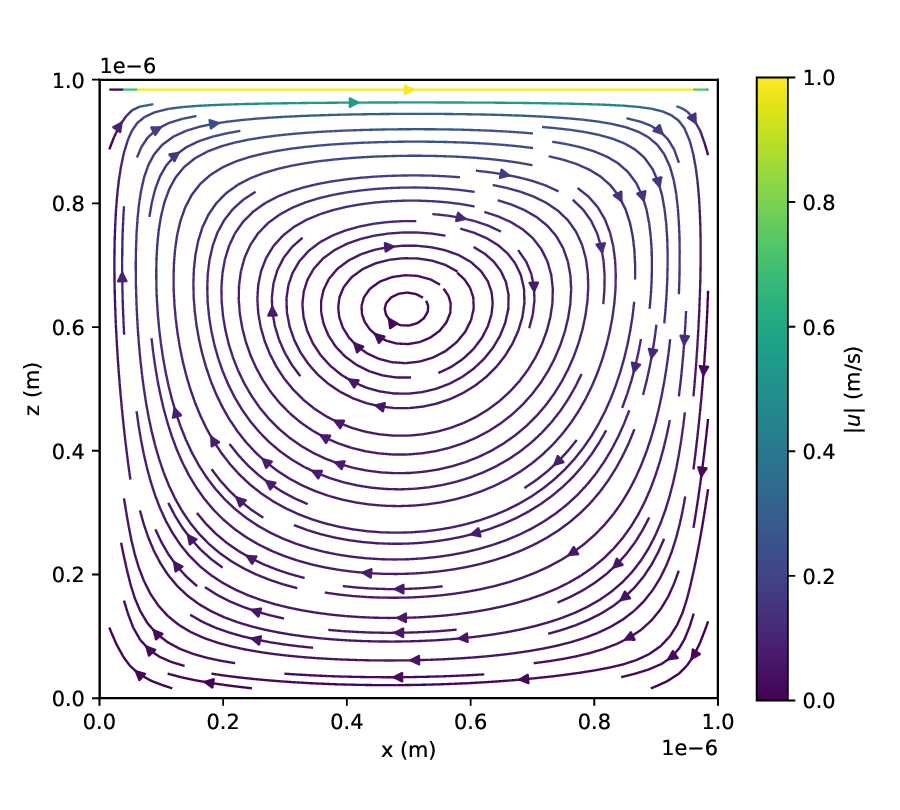}
	\includegraphics[keepaspectratio=true, angle=0, width=0.48\textwidth]{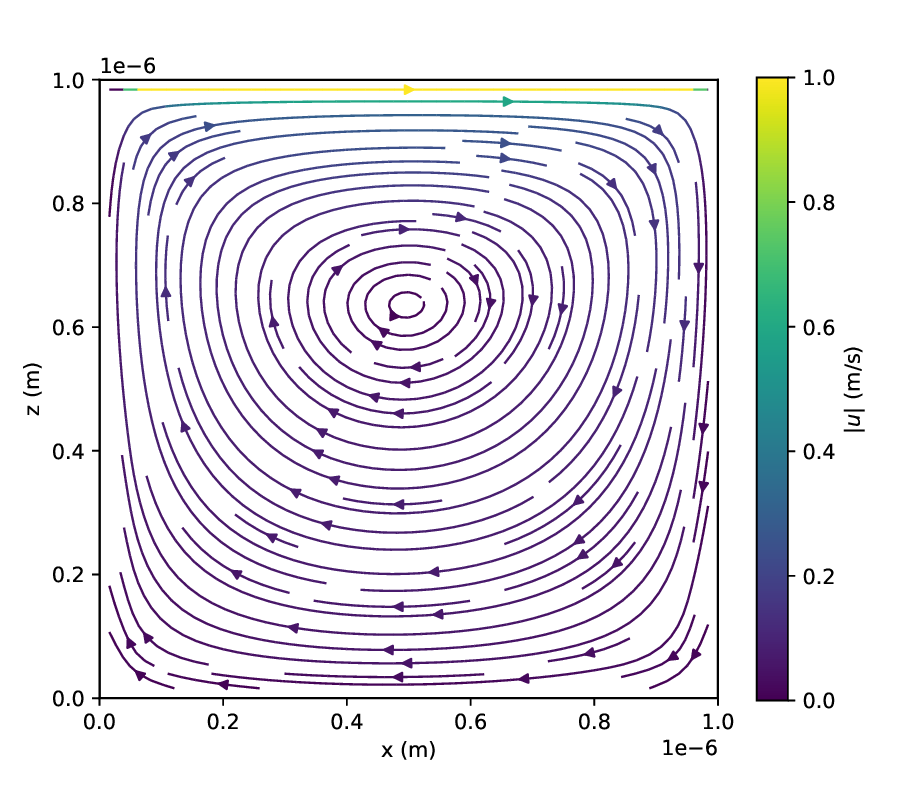}
	\caption{Three dimensional velocity stream plot in $xz$ plane at half of the $y$-axis for Knudsen numbers and Reynolds numbers $Kn = 1.1, Re = 2.968 \times 10^{-4}$ (left) and $Kn = 0.11, Re = 2.968 \times 10^{-5}$ (right).}
	\label{xzplane}
	\centering
\end{figure}

In this case we consider a cube  $[0, L] \times [0, L] \times [0, L] $ with  $L=1\times 10^{-6}$.  The cloud is generated as $N$ regular grid points, an arrangement that is lost after a few time steps as
the points move with the flow.
The initial parameters of the Maxwellian are the same as in the previous 2d case. We use diffuse reflection boundary conditions again with fixed temperature  $T_0$  on all boundaries. On the top wall $z=L$, we prescribe a non-zero velocity in the $x$-direction given by	$\bf {U} =  (1,0,0).$
On the other five walls, all velocities are zero. In Fig. \ref{3dcavity} we have plotted the velocity vectors for about $27,000$ particles for Knudsen numbers $Kn = 1.1$ and $0.11$. Moreover, in Fig. \ref{xzplane} we have plotted the velocity vectors in the $xz$ plane at $y=0.5\times 10^{-6}$ for both Knudsen numbers. One observes that the results are qualitatively the same for both Knudsen numbers and similar to the two-dimensional case.	
With $N_v=20$ the single distribution $f$ occupies $N_v^3 = 8000$ velocity nodes per particle, i.e.\ $128$\,KB. The $100^3$ grid requires approximately 120\,GB (see Table~\ref{tab:slablimits}). The 3D transport evolves the single distribution $f$ on the full $N_v^3$ velocity grid, as described in Section~\ref{sec:cuda_upwind}.  Each velocity node carries its resolved $v_3$, so the $z$-transport follows the true characteristic in all three directions with no bulk-velocity closure.  The per-step cost is therefore $O(N_v^3)$, which sets the memory footprint reported above. With the velocity cutoff set to $C_{\rm cut}=5.5$, the initial Maxwellian is captured with a quadrature error of less than $10^{-3}$\, K on the temperature.

To compare the speedup, we again consider several resolutions
(Table~\ref{tab:bench3dcav}) at fixed velocity resolution. For one million
grid points with $N_v = 16$, the GPU is $1087$ times faster than the serial
CPU, $209$ times faster than the workstation and $79$ times faster than the
HPC node. For $N = 80^{3}$ with $N_v = 20$, a complete run at the uniform
benchmark step would take over eleven years in serial, two and a half years
on the workstation and some $326$ days on the HPC node, against roughly four
days on an A100. At the recommended step of Section~\ref{sec:timestep} the
same run requires $1.2$ hours on the GPU, against $12.6$ days and $4.5$ days
on the workstation and HPC respectively.

\begin{table}[!ht]
	\begin{center}
		\begin{tabular}{|r|r|r|r|c|c|c|c|c|c|}
			\hline
			\multirow{2}{*}{$N_x$} & \multirow{2}{*}{$N$} & \multirow{2}{*}{$N_v$} & \multicolumn{4}{c|}{Time per step (s)} & GPU full run \\ \cline{4-7}
			& & & Serial & Workstation (P=8) & HPC (P=8) & GPU  & (h) \\ \hline
			64 & 262,144 & 16 & 1926.34  & 417.28 & 131.79  & 1.73 & 24.03 \\     
			80 & 512,000 & 16 & 3665.09  & 720.94&  264.31 & 3.31 & 46.1 \\
			100 & 1,000,000 & 16 & 7011.26  & 1345.84 & 508.43  & 6.45 & 89.7 \\                           
			\hline
			64 & 262,144 & 20 & 3237.75  & 800.48 & 279.14   & 3.4 & 47.3 \\     
			80 & 512,000 & 20 & 7039.61  & 1591.34 & 563.43  & 6.5 & 90.4 \\
			\hline
		\end{tabular}
		\caption{3D driven cavity: measured time per step for the serial, MPI
			and A100 GPU on Google Colab Pro. ($\Delta t=10^{-13}, t_{\rm final}=5\times 10^{-9}$, 50,000 steps).}
		\label{tab:bench3dcav}
	\end{center}
\end{table}

\subsection{Moving rigid bodies}\label{sec:rigidbodies}

A rigid body immersed in the gas is represented by a layer of surface
particles. Each carries the outward body normal ${\bf n}$, the local wall
velocity ${\bf U}_w = {\bf V} + \boldsymbol{\omega}\times{\bf r}$, with
${\bf r}$ the position relative to the centre of mass, and the wall
temperature, and participates in the diffuse-reflection boundary condition
exactly as the container walls do \cite{Klaas}. The hydrodynamic load on the body is the
surface integral of the kinetic stress evaluated relative to the wall
velocity,
\begin{eqnarray}
	\boldsymbol\Phi_k = \int ({\bf v} - {\bf U}_w)\otimes({\bf v} - {\bf U}_w) f\, d{\bf v}, \quad {
		\bf F} = -\sum_k A_k\, \boldsymbol\Phi_k\,{\bf n}_k, \quad {\bf M} = \sum_{k} A_k\, {\bf r}_k \times (-\boldsymbol\Phi_k\,{\bf n}_k) 
\end{eqnarray}
where the sum runs over the surface particles with area
weights $A_k$ and lever arms ${\bf r}_k$ from the centre of mass. The
Newton--Euler equations are advanced with a Heun (predictor-corrector)
update of the interface loads: a first load evaluation at the current wall
state predicts the body state, the loads are re-evaluated with the
predicted wall velocities on the same frozen distribution, which captures
the drag and torque Jacobian $\partial {\bf F}/\partial {\bf V}$ of the
diffuse wall, and the state is advanced with the averaged load,
\begin{eqnarray}
	{\bf V}^{n+1} &=& {\bf V}^{n} +
	\frac{\Delta t}{m}\,\frac{{\bf F}({\bf V}^n) + {\bf F}({\bf V}^{*})}{2},
	\\
	\boldsymbol{\omega}^{n+1} &=& \boldsymbol{\omega}^{n} +
	\frac{\Delta t}{I}\,\frac{{\bf T}(\boldsymbol{\omega}^n) + {\bf
			T}(\boldsymbol{\omega}^{*})}{2}
\end{eqnarray}
A staggered update with a load lagged by one full step
oscillates around the correct equilibrium with growing amplitude for a
freely rotating body. The Heun coupling removes the lag without any
relaxation parameter and is exact at load equilibrium.

\subsubsection{Two-dimensional driven cavity with a moving body}
\label{sec:cavity2d_body}

The two-dimensional cavity occupies $[0,L]^{2}$ with $L = 1 \times 10^{-6}\,m$ filled
with argon at $270$\,K and bounded by diffuse walls, the lid at $y = L$ being
driven with a parabolic profile of peak speed $10$\,m\,s$^{-1}$. A square of
side $0.075L$, of density ten times that of the gas, is released at rest with
its centre at $(0.6,\,0.7)L$ and its edges aligned with the axes. The centre
is displaced from the cavity centre deliberately, on the axis of the primary
vortex; the net torque very nearly vanishes, so a body released there would neither spin nor migrate
\begin{figure}[!ht]
	\centering
	\includegraphics[keepaspectratio=true, angle=0, width=\textwidth]{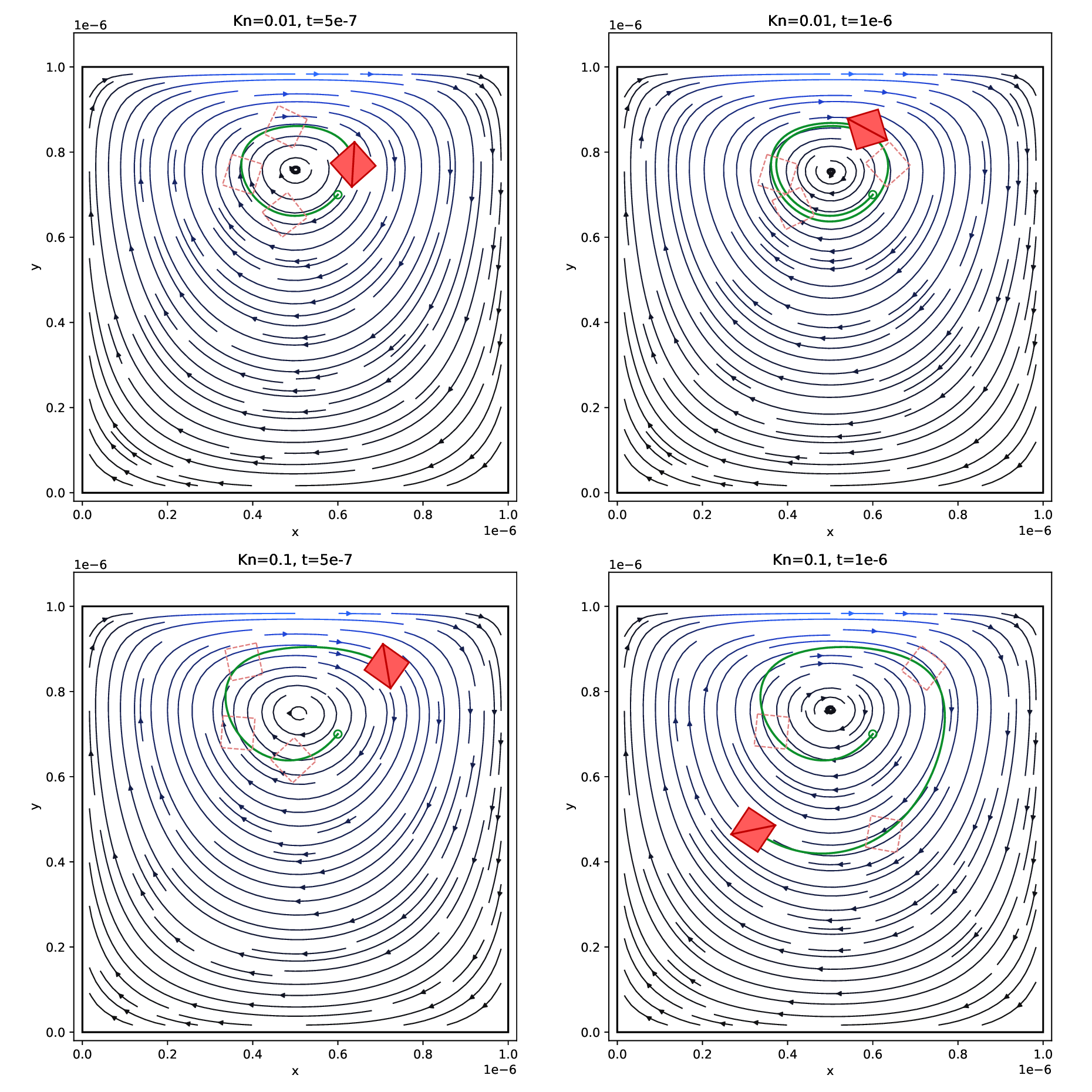}
	\caption{Freely moving square in the driven cavity: streamlines coloured
		by speed, the body-centre trajectory travelled so far (green), earlier
		body poses (dashed) and the body at its instantaneous position and
		orientation (red), for $Kn = 0.01$ (top) and $Kn = 0.1$ (bottom) at
		$t = 5\times10^{-7}$ (left) and $t = 10^{-6}$ (right).}
	\label{cavitybody}
\end{figure}

The square is set into rotation by the torque of the recirculating flow and
is simultaneously advected around the primary vortex. Over the simulated $10^{-6}$\,s it
traverses a path of $1.374L$ while its net displacement is only $0.152L$. The body orbits the vortex and returns close to its point of release rather than drifting to a wall. 

The rotation does not settle. The body accumulates $616^{\circ}$ over the run
and the angular velocity grows from a
first-quarter mean of $7.4\times10^{6}$\, rad\,s$^{-1}$ to a last-quarter mean
of $1.4\times10^{7}$\,rad\,s$^{-1}$, with no plateau. A body held at a fixed station would spin up towards the local shear rate and stop there, whereas one carried around the vortex meets a different shear at every station and is continually re-accelerated.
Figure~\ref{cavitybody} shows the motion for both Knudsen numbers, with the
streamlines of the driving vortex, the centre trajectory traced up to each
panel time, and the body at its instantaneous position and orientation,
faded outlines marking earlier poses.

\subsubsection{Three-dimensional driven cavity with a moving body}
\label{sec:cavity3d_body}

The three-dimensional case occupies $[0,L]^{3}$ under the same wall
temperature and lid profile, with the lid again driving the $y = L$ face. A cube
of side $0.15L$ and density ten times that of the gas is released at rest,
axis-aligned, with its centre at $(0.6,\,0.7,\,0.5)L$ -- offset within the
primary vortex as in two dimensions, and at mid-span in $z$. The gas is lighter than
in the two-dimensional case, at a density giving $\lambda/\Delta x \approx
0.96$, so that the Knudsen layer is resolved on the grids used here; the two
cases therefore differ in Knudsen number as well as in dimension.

\begin{figure}[!ht]
	\centering
	\includegraphics[keepaspectratio=true, angle=0, width=\textwidth]{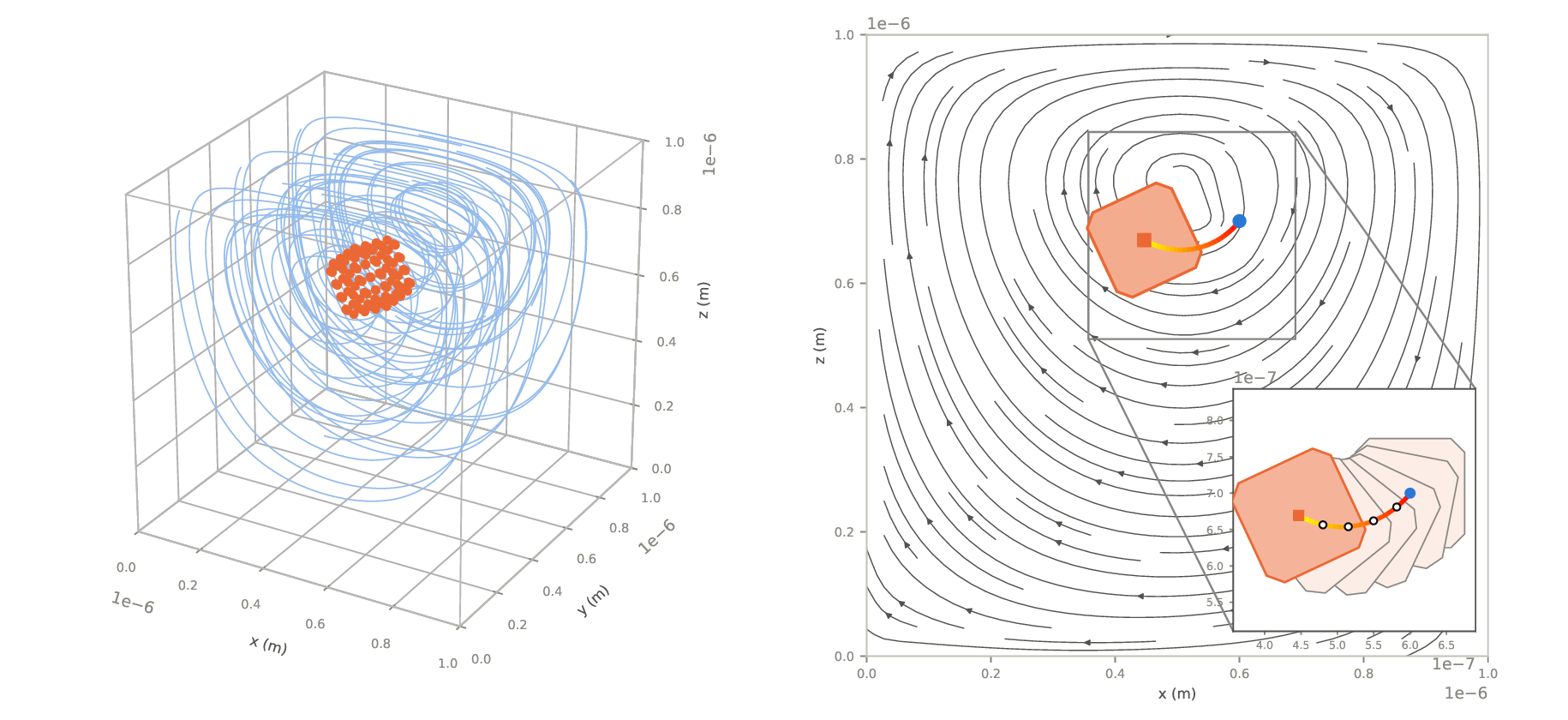}
	\caption{Free cube in the 3-D driven cavity, $N_x = 40$, $N_v = 10$,
		$t = 1.8\times10^{-7}$\,s. (a) Three-dimensional streamlines of the
		carrying flow; arrows on the top face mark the driven lid, which slides
		in $+x$; the cube is drawn from its surface particles at the final pose.
		(b) The $y = L/2$ slice: in-plane streamlines, the body-centre path from
		release (circle) to final position (square), and the cube's cross-section.
		The inside figure magnifies the travelled region.}
	\label{cavity3dbody}
\end{figure}	

\begin{figure}
	\centering
	\includegraphics[width=0.55\textwidth]{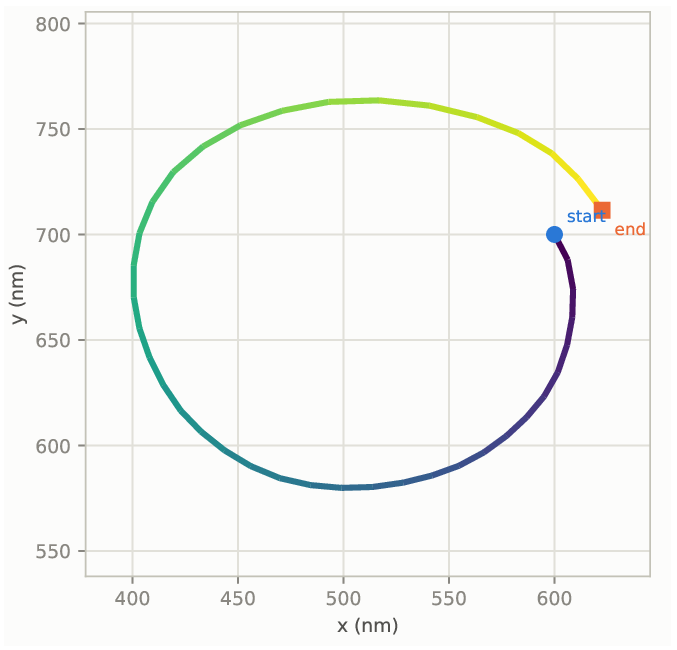}
	\caption{Body-centre trajectory over a complete revolution, $N_x = 40$,
		$N_v = 10$. The cube is released at the circle and returns to the square
		after $t \approx 1.17\times10^{-6}$\,s, by which time it has turned
		through $360^{\circ}$.}
	\label{fullrotpath}
\end{figure}

The motion is very nearly planar. Over the integrated interval, the cube's
spanwise coordinate stays within $0.003L$ of mid-span while its in-plane
centre moves by $0.113L$, so the body is carried by the primary vortex in the
$x$--$y$ plane and is not measurably displaced along the lid-parallel
direction. This is consistent with the driving vortex being essentially
two-dimensional away from the end walls. Figure~\ref{cavity3dbody} shows the resulting motion: panel~(a) the
three-dimensional streamlines of the carrying flow, whose axis lies normal to
the direction of lid motion, with the cube at its final pose; panel~(b) the
mid-span slice, carrying the body-centre path and the cube's cross-section.
The magnified inset resolves the travelled region, where the successive
cross-sections make the rotation plain while the centre advances only a
fraction of the body's own width. Carried far enough, the in-plane motion closes over the longer interval of
Figure~\ref{fullrotpath}, integrated on a coarser grid, the centre completes
a circuit of the primary vortex and returns to within $0.02L$ of its release
point, having turned through a full revolution.

Where the two cases differ sharply is in the approach to a steady spin. At a
matched instant of $t = 1.8\times10^{-7}$\,s the square has turned
$74.2^{\circ}$ with its angular velocity still climbing steeply (a factor
of $2.8$ across the interval) whereas the cube has turned $58.4^{\circ}$ at
a rate that is already constant to within $6\%$, its first- and
last-quarter means differing by a factor of $0.94$. The cube thus reaches a
near-steady rotation early and holds it, at $5.5\times10^{6}$\,rad\,s$^{-1}$,
while the square is still accelerating. At this rate, a complete revolution
requires $t \approx 1.13\times10^{-6}$\,s.

The accumulated rotation is the natural scalar against which to test grid
convergence, since it integrates the surface-stress torque over the whole
trajectory and is insensitive to the phase of the orbit. A four-grid study at
$N_x = 30,\,40,\,50,\,60$ places the converged value at approximately
$65.3^{\circ}$ at $t = 1.8\times10^{-7}$\,s, with $N_x = 40$ some $11\%$ below
it and $N_x = 50$ within $2.7\%$. The result is, by contrast, almost
independent of the velocity grid: refining from $N_v = 10$ to $N_v = 20$ at
$N_x = 40$ changes the accumulated rotation from $58.34^{\circ}$ to
$58.42^{\circ}$, a difference of $0.14\%$. The body's motion is therefore
controlled by the spatial resolution of the surface-stress quadrature and of
the Knudsen layer around it, not by the resolution of the velocity space.

\begin{table}[H]
	\centering
	\begin{tabular}{|c|c|c|c|c|}
		\hline
		$N_x$ & $N$ & $N_v$ & GPU & HPC  \\
		\hline
		$80$ & 512,000 & $20$ & 37.16& 614.58\\
		$100$ & 1,000,000 & 16&35.76 & 569.8\\
		\hline
	\end{tabular}
	\caption{3D lid-driven cavity with a free cubic body: measure time per step for A100 GPU and a single HPC rank with $20$ OpenMP threads,
		each averaged over $100$ time steps.}
	\label{tab:3drivecavfreebody}
\end{table}

Table~\ref{tab:3drivecavfreebody} reports the measured cost per time step
with the body present, for the two largest configurations reached here. The
GPU completes a step about sixteen times faster than the HPC node in both
cases, against the seventy- to ninety-fold advantage measured at the same
resolutions without a body (Table~\ref{tab:bench3dcav}). The device time per
step rises by a factor of $5.5$ to $5.7$ once the body is embedded while the
kinetic work is unchanged, which places the host--device transfer of the
point cloud, rather than the kinetic update, as the cost the moving-body path
adds.

\subsubsection{Computational cost}
Table~\ref{tab:cubeprofile} breaks one time step of the three-dimensional
cavity with a free cubic body ($40^{3}$ particles, $N_v = 20$, that is, $63,880$
particles once the body is embedded) into its constituent costs. From the table, it is evident that the transport kernel dominates the arithmetic cost. A moving
body gives the step two cost centres rather than one: the kinetic update on
the device, and an Arbitrary Lagrangian--Eulerian pass on the host that
maintains the point cloud as the body sweeps through it. On this hardware the
host pass is the larger of the two by a wide margin, at $81.1\%$ of the step
against $17.5\%$ for the GPU kernels; the two sum to $4.32$\,s against a
measured $4.39$\,s, a residual of $1.4\%$. The distribution arrays occupy
$7.65$\,GB, comfortably within the $40$\,GB device.

\begin{table}[!ht]
	\begin{center}
		\begin{tabular}{|l|r|r|}
			\hline
			\textbf{Component} & \textbf{Time per step (ms)} & \textbf{Share (\%)} \\ \hline
			\multicolumn{3}{|c|}{\textbf{GPU kernels}} \\ \hline
			Transport (3D, $f$ on $N_v^3$)             &    642.88 &  14.66 \\
			Commit ($g_t \rightarrow g$)               &     34.11 &   0.78 \\
			Mass correction                            &     19.49 &   0.44 \\
			Collision (implicit BGK)                   &     16.16 &   0.37 \\
			Wall boundary condition                    &     15.27 &   0.35 \\
			Admissibility (NAD) clamp                  &     11.61 &   0.26 \\
			Transport-buffer prefill                   &      6.09 &   0.14 \\
			Mass reduction                             &      5.36 &   0.12 \\
			Positivity fix                             &      5.15 &   0.12 \\
			Moments (post-collision)                   &      4.76 &   0.11 \\
			Moments calculation                        &      4.74 &   0.11 \\ \hline
			Subtotal                                   &    \textbf{765.61} &  \textbf{17.46} \\ \hline
			\multicolumn{3}{|c|}{\textbf{Host ALE pass (moving body)}} \\ \hline
			Point cloud download (device $\rightarrow$ host) &   1739.30 &  39.65 \\
			Point cloud upload (host $\rightarrow$ device)   &   1728.50 &  39.41 \\
			Neighbour-list rebuild                     &     57.20 &   1.30 \\
			Merge close particles                      &     14.00 &   0.32 \\
			Remove overrun particles                   &     10.80 &   0.25 \\
			Rigid-body update                          &      7.60 &   0.17 \\
			Refill vacated voxels                      &      0.20 &   0.00 \\ \hline
			Subtotal                                   &   \textbf{3557.60} &  \textbf{81.11} \\ \hline
			Measured step                              &   4386.18 & 100.00 \\
			\hline
		\end{tabular}
		\caption{Kernel and host profile of one time step, 3D lid-driven
			cavity with a free cubic body, $40^{3}$ particles, $N_v = 20$. GPU times are means over fifteen steady-state steps; host times are means over one hundred steps; the measured step is the mean over the $110$-step run.}
		\label{tab:cubeprofile}
	\end{center}
\end{table}

Within the device kernels, most of the cost is for the semi-Lagrangian
transport of the single distribution $f$, at $84.0\%$ of GPU time. Every
other kernel costs under $35$\, ms, including the implicit relaxation, which
touches the same $N \cdot N_v^{3}$ nodes for $0.4\%$ of the step but performs
no neighbour gathers. The contrast isolates the neighbour-gathering
reconstruction, rather than the node count, as what the kinetic update
actually pays for, and motivates the register-resident, compile-time-sized
least squares solve of Section~\ref{sec:cuda_upwind}.

The host pass, however, is not dominated by the remeshing it exists to
perform. Uploading and downloading the cloud together are $79.1\%$ of the
step, while removing overrun particles, merging close pairs, refilling
vacated voxels and integrating the body amount to $0.7\%$ combined, and the
neighbour-list rebuild to a further $1.3\%$. The moving-body treatment cost is therefore almost negligible.

This balance is a consequence of the device being fast at double precision
rather than a fixed property of the scheme. The kinetic update is served by
full-rate FP64 and falls to a sixth of the step, while the transfers are
bounded by host-device bandwidth and are untouched by it; on hardware with
reduced-rate FP64 the same decomposition would instead be dominated by the
transport kernel. On a machine of the present class the transfer schedule is
thus the binding constraint, and restricting the round trip to the steps on
which the cloud actually changes, rather than every step, is the single
highest-value optimisation available to the moving-body path.

\section{Conclusion and Outlook}
\label{sec:conclusion}
In this paper, we have presented an Arbitrary Lagrangian-Eulerian method for
the simulation of the BGK equation with GPU parallelization. Besides the ALE
approach, the method is based on a semi-Lagrangian transport step with least
squares reconstruction on the moving point cloud. Two- and three-dimensional
driven cavity test cases are performed, together with
freely moving rigid bodies. In the two-dimensional driven cavity, the GPU gives a speedup of about $757$ times over a sequential CPU computation in two dimensions and
about $1087$ times in three dimensions. Compared with a distributed-memory CPU run, it remains about $50$ times faster than 16 MPI ranks and $29$ times faster than $40$ ranks in two dimensions, and about $209$ and $79$ times faster respectively in three. With a freely moving body, the per-step host--device transfer required by the
particle management reduces the advantage to about sixteen times over the HPC
node, which remains decisive but leaves the transfer schedule, rather than
the kinetic update, as the binding constraint.

The meshfree transport step has been implemented as an efficient
node-parallel CUDA kernel for both 2D and 3D physical space. The key
algorithmic ingredients are: (1) parallelism over particle-velocity-node
pairs rather than over particles alone, (2) register-resident least squares
reconstruction with compile-time system sizes, using a pivot-free Cholesky
solve of the symmetric positive definite normal matrix in 3D, and (3)
atomic-free block reductions for the global conservation correction. In
three dimensions, no velocity reduction applies, so a single distribution
$f$ is evolved on the full $N_v^3$ grid with the $z$-transport carried by
the resolved velocity node. Memory management relies on a
structure-of-arrays device cloud whose footprint is exactly twice the
distribution storage, and particle management runs on the host, so no
device-side scratch copies are required. Experiments indicate that the velocity cutoff must satisfy
$C_{\rm cut}\geq 5$ thermal speeds, which keeps the discrete quadrature of
the initial Maxwellian within $10^{-3}$\,K of the prescribed wall
temperature. Smaller cutoffs introduce a systematic mismatch that seeds a
slowly growing instability.

In future work, the scheme will be extended to higher-order approximations
using MPI and a multi-GPU optimized code. Furthermore, the case of
gas mixtures \cite{GR} will be considered, and the rigid-body coupling
presented here extended to several interacting bodies and to non-isotropic
inertia tensors.

\subsection*{Data and Code Availability}
The complete code is available in the GitHub repository. 
\url{https://github.com/panchiittp/ALEBGK4Paper/}. Also, the Google Colaboratory code is available to copy, run and test the manuscript: \url{https://colab.research.google.com/drive/1qXJ184S93gNI-kWk00_RpKJYvyGCqAsk?authuser=3#scrollTo=c9}

\subsection*{Declaration of generative AI and AI-assisted technologies in the
	writing process}

During the preparation of this work, the authors used Grammarly for grammar
and typographical corrections, and Claude Code Pro (Anthropic) for the following
purposes: preparing the shell and PBS job scripts used to launch the
solver; configuring and executing runs on the authors' hardware unattended;
and generating explanatory comments for existing source code. The numerical
method, the simulations, the results and all scientific conclusions are the
authors' own. After using these tools, the authors reviewed and edited the content as needed and take
full responsibility for the content of the publication.

\subsection*{Acknowledgments}
This  work is supported
by the DFG (German Research Foundation) under Grant No. KL 1105/34 and DAAD-DST India 2022-2024, Project-ID: 57622011 \& DST/INT/DAAD/P-10/2022 and the European Union's Framework Program for Research and Innovation HORIZON-MSCA-2021-DN-01 under the Marie Sklodowska-Curie Grant Agreement Project 101072546 DATAHYKING.

\bibliographystyle{plain}
\bibliography{References}

	
\end{document}